\documentclass[11pt,english]{article}

\let\origrmdefault\rmdefault
\usepackage[math]{iwona}
\renewcommand{\rmdefault}{\origrmdefault}
\usepackage[T1]{fontenc}
\usepackage[latin9]{inputenc}
\usepackage{babel}
\usepackage{array}
\usepackage{amsmath}
\usepackage{amsthm}
\usepackage{amssymb}
\usepackage{stmaryrd}
\usepackage{graphicx}
\usepackage{geometry}
\usepackage{setspace}
\usepackage[authoryear]{natbib}
\usepackage[pdfusetitle,
 bookmarks=false,
 breaklinks=false,pdfborder={0 0 1},backref=false,colorlinks=false]
 {hyperref}

\makeatletter

\providecommand{\tabularnewline}{\\}

\newcommand{\lyxaddress}[1]{
	\par {\raggedright #1
	\vspace{1.4em}
	\noindent\par}
}

\@ifundefined{date}{}{\date{}}
\usepackage{babel}

\usepackage{pdflscape}
\usepackage{amsthm}
\newtheorem{theorem}{Theorem}\newtheorem{proposition}[theorem]{Proposition}\newtheorem{remark}[theorem]{Remark}

\makeatother

\begin{document}
\title{Quadratic invariants and Hamiltonian structure in coupled gyrostat
low-order model hierarchies}
\author{Ashwin K Seshadri\textsuperscript{1} and S Lakshmivarahan\textsuperscript{2}}
\maketitle

\lyxaddress{\textsuperscript{1}Centre for Atmospheric and Oceanic Sciences and
Divecha Centre for Climate Change, Indian Institute of Science, Bangalore
560012, India. Email: ashwins@iisc.ac.in.}

\lyxaddress{\textsuperscript{2}Emeritus faculty at the School of Computer Science,
University of Oklahoma, Norman, OK 73012, USA. Email: varahan@ou.edu.}

\subsection*{Declarations of interest: none}

\pagebreak{} 
\begin{abstract}
Coupled gyrostat low-order models (GLOMs) are energy-conserving cores
of Galerkin-truncated fluid and geophysical systems, including Rayleigh--B�nard
convection and vorticity dynamics. A single gyrostat always possesses
two quadratic invariants; when gyrostats are coupled, the number and
geometry of invariants vary sensitively with model configuration,
influencing the effective dimension of the dynamics, nonlinear stability,
and statistical equilibria. We provide a systematic theory of this
dependence. For sparse nested hierarchies of $K$ gyrostats ($M=2K+1$
modes, no linear feedback), the number of independent quadratic invariants
is exactly $(M+1)/2$; for general GLOMs with all parameters nonzero,
energy is the only guaranteed invariant.

The standard algebraic approach to finding invariants does not scale
with model size. We show instead that many GLOMs admit a non-canonical
Hamiltonian structure, with quadratic invariants recoverable as Casimir
functions of an explicitly constructible Poisson matrix. The Hamiltonian
structure imposes precise, computationally verifiable constraints
on the nonlinear coefficients. For Hamiltonian hierarchies, Casimir
gradients project consistently across models of increasing complexity,
so that invariants are compatible under restriction to subspaces.
The clear interpretation of these models thus enables consistent application
of Hamiltonian dynamics across low-order model hierarchies. 
\end{abstract}

\section{Introduction}

The Lorenz system \citet{Lorenz1963}, derived by Galerkin truncation
of the Boussinesq fluid equations, is an early and widespread example
of deterministic chaos in a low-order model. Its conservative core,
obtained by removing forcing and dissipation, is a three-dimensional
gyrostat that conserves two independent quadratic invariants: energy
and a second quantity that together confine every trajectory to a
closed curve. As a result, chaos cannot occur in this model. Adding
a second gyrostat to obtain a four-dimensional coupled system changes
this picture entirely. Depending only on which linear feedback terms
are present or absent in the model, the coupled system can have either
one, two, or three invariants, as we shall demonstrate. In general
only the single invariant of energy is guaranteed for coupled gyrostat
models. With only this single invariant, the dynamics is three-dimensional
on the energy surface and irregular motion is possible, while for
the case of three invariants the dynamics collapses onto a curve.
This sensitivity of the invariant number and their structure, and
therefore of dynamical complexity, to the configuration of the model
in such systems is the central phenomenon explained by this paper.

Low-order models (LOMs) obtained by Galerkin projection of partial
differential equations (PDEs) are standard tools in geophysical fluid
dynamics and climate science. When external forcing and dissipation
are stripped away from these models, these LOMs expose a conservative
core whose behaviour reflects the intrinsic nonlinear dynamics of
the parent PDE. Remarkably, many such conservative LOMs, arising from
systems as wide-ranging as the barotropic vorticity equation, the
Boussinesq equations, and quasi-geostrophic potential vorticity equation,
and many others, can be written as systems of coupled gyrostats, where
each gyrostatic component preserves energy and generates a skew-symmetric
linear contribution to the vector field. This representation, developed
systematically by \citet{Oboukhov1975}, \citet{Gluhovsky1997}, and
\citet{Gluhovsky1999}, encompasses the Lorenz system \citep{Lorenz1963},
barotropic models \citep{Charney1979}, wave--mean flow interactions
\citep{Swart1988}, and turbulent convection \citep{Howard1986},
among others. We refer to such systems collectively as Gyrostat Low-Order
Models (GLOMs).

A single gyrostat always conserves two independent quadratic invariants:
energy and a second quantity whose explicit form depends on the gyrostat
subclass \citep{Gluhovsky2006,Seshadri2023}. This second invariant
is physically analogous to the squared angular momentum of the physical
gyroscope and ensures that the dynamics of the single gyrostat is
inevitably periodic in the absence of any forcing or dissipation.
When $K$ gyrostats are coupled to form a GLOM with $M\geq4$ modes,
the number of invariants can range from one (energy alone) to $(M+1)/2$
for sparse coupling with odd-valued $M$ as shown below. Generally
the number of invariants depends sensitively on coupling configuration
and the presence of linear feedback terms.

The implications of this wide range of the number of invariants extend
well beyond the question of whether individual trajectories can be
irregular or must necessarily be quasiperiodic or periodic. Each independent
invariant $C_{j}$ constrains every trajectory to the level set $\{C_{j}(\mathbf{x})=c_{j}\}$,
and together the $n_{c}$ invariants confine the dynamics to a smooth
submanifold $\mathcal{L}_{\mathbf{c}}\subset\mathbb{R}^{M}$ of dimension
$M-n_{c}$. For a non-canonical Hamiltonian system, this submanifold
forms a symplectic leaf of the Poisson structure \citep{Marsden1999}
and, although the Poisson bracket is degenerate on $\mathbb{R}^{M}$
it becomes nondegenerate when restricted to each leaf, and the dynamics
on the leaf is thereby a classical Hamiltonian system in case $M-n_{c}$
is even. Then the full apparatus of classical Hamiltonian dynamics,
including KAM theory \citep{Ott1993}, Birkhoff normal forms, and
related tools for nearly integrable systems apply to the restricted
dynamics but not to the ambient dynamics in $\mathbb{R}^{M}$.

The symplectic leaf structure that accompanies the non-canonical Hamiltonian
reduction has additional effects that are directly relevant to the
geophysical applications motivating the systematic treatment of GLOMs.
First, it enables nonlinear stability analysis via Arnold's method,
wherein a steady state can be proven to be nonlinearly stable whenever
the second variation of the constrained Hamiltonian $H+\sum_{j}\lambda_{j}C_{j}$
is positive definite on perturbations tangent to the Casimir leaf,
a condition that requires the existence of sufficiently many Casimirs
that eliminate neutrally stable directions. This method has been applied
systematically to geophysical fluid equilibria by \citep{Holm1985}
and \citep{Shepherd1990}; for GLOMs it implies that the presence
of more invariants can in fact support a richer nonlinear stability
theory. Second, even when dynamics on $\mathcal{L}_{\mathbf{c}}$
is chaotic, the long-time statistical state is constrained to distributions
that are supported on $\mathcal{L}_{\mathbf{c}}$. An important consequence
is that the accessible statistical equilibria are determined by which
invariants exist, not by the detailed structure of trajectories within
the manifold. This is a key mechanism underlying the observation \citep{Kraichnan1967}
that enstrophy conservation reverses the direction of energy cascade
in two-dimensional turbulence, and is formalised in its equilibrium
statistical mechanics \citep{Miller1990,Robert1991}. Third, when
forcing and dissipation are added, different invariants decay at different
rates, producing the phenomenon of selective decay \citep{Matthaeus1980}
toward specific invariant manifolds. Then the the Casimir structure
of the conservative core can shape the statistics of the forced-dissipative
system even though Casimirs are no longer conserved in the presence
of dissipation.

A further consequence concerns consistency across various levels of
truncation. A Galerkin truncation of a PDE does not in general preserve
the Poisson structure of the parent equation \citep{Zeitlin1991},
and a LOM whose conservative core has different invariant structure
from the full model will predict rather different (and possibly incompatible)
statistical equilibria and stability properties. Hamiltonian GLOM
hierarchies, the main constructive output of this paper, can be demonstrated
to have truncation families for which invariant structure is consistent
across models in the hierarchy: the Casimirs of the $K$-gyrostat
model are the restrictions of the Casimirs of the $(K+1)$-gyrostat
model, so the symplectic geometry, the nonlinear stability theory,
and the accessible statistical equilibria can be expected to be similarly
consistent throughout the hierarchy.

Identifying and counting the invariants of GLOMs, and determining
how non-canonical Hamiltonian structure enables this consistency,
is the central analytical problem we address in this paper.

Invariants are conventionally found by requiring $\dot{C}(\mathbf{x})=0$,
which for polynomial systems yields a linear system in the coefficients
of $C$. This approach is tractable for individual low-dimensional
cases but does not scale: the number of parameter subclasses grows
as $2^{6K}$ for $K$ gyrostats, even once the coupling between gyrostats
and modes is fixed, since each gyrostat has $6$ parameters. Identifying
which mixed quadratic terms vanish requires case-specific analysis
for each configuration, and detecting algebraic dependence among the
various candidate invariants adds further difficulty to a general
analysis of this problem. Symmetries of the vector field can often
reduce the calculation (also see \citet{Seshadri2023}), but are present
only in subclasses where several linear feedback parameters simultaneously
vanish, which are precisely the special cases that are already tractable
on the conventional approach. Therefore for the generic GLOM, the
conventional approach provides no practical or computational simplifications.

A more tractable framework is provided by non-canonical Hamiltonian
mechanics. When a GLOM admits a Poisson structure, involving a skew-symmetric
matrix $\mathrm{J}(\mathbf{x})$ that satisfies the Jacobi identity
and generates the vector field via $\dot{x}_{i}=\mathrm{J}_{ij}\partial H/\partial x_{j}$,
new possibilities for analysis of consistency between invariants arise.
This is because invariants appear as Casimir functions from the nullspace
of $\mathrm{J}$. \citet{Arnold1969} first showed that the Euler
gyrostat has this structure; Gluhovsky \citep{Gluhovsky2006} extended
it to the Volterra gyrostat. In the fluid-dynamical context, Casimirs
are analogues of conserved integral quantities such as enstrophy or
potential vorticity \citep{Morrison1996}: they are conserved regardless
of initial condition and define invariant submanifolds of the finite-dimensional
phase space of the low-order model. Whether a GLOM possesses Casimirs,
and how many, therefore has direct physical consequences for the dynamics
of the truncated model and its various properties described above.
A systematic theory of which GLOMs admit non-canonical Hamiltonian
structure, and how this structure governs the invariant count, has
not previously been developed.

In this paper we attempt to provide that theory. We show that the
Jacobi condition on $\mathrm{J}$, expressed as a sum of gyrostat
contributions, yields explicit parameter constraints for Hamiltonian
structure that can be verified computationally and that propagate
through nested model hierarchies via a recurrence (Proposition~\ref{prop:recurse}).
In the general case with all parameters nonzero, energy is the only
invariant that is assured regardless of the coupling configuration
between gyrostats and modes. For Hamiltonian hierarchies of both nested
and fully coupled type, including models of 2D and 3D Rayleigh--B�nard
convection, we show that Casimir gradients project consistently under
restriction to subspaces, so that invariants are evidently compatible
across models of increasing complexity (Theorem~\ref{prop:consist}).
This consistency property provides a geometric criterion for building
physically faithful reduced-order model hierarchies where, additionally,
the various tools of Hamiltonian dynamics can be deployed.

\subsection*{Main results}

The central finding of this paper is that non-canonical Hamiltonian
structure provides both a diagnostic and a constructive tool for the
analysis of invariants of GLOMs, which circumvents the combinatorial
intractability of the standard algebraic approach to finding invariants
for general non-Hamiltonian systems.

We first show that the number of quadratic invariants in GLOMs ranges
widely, from one (energy alone) to $(M+1)/2$ for special configurations
and with parameter restrictions (where $M$ is the number of modes),
and that this range depends sensitively on the coupling configuration
and especially on whether and which linear feedback terms are present.
Making additional parameters nonzero cannot increase the invariant
count (Proposition~1), so the minimum number of invariants occurs
generically while the maximum is achieved only in special parameter
regimes. For sparse nested hierarchies with $M=2K+1$ modes and no
linear feedback, the number of independent quadratic invariants is
exactly $(M+1)/2$ (Theorem~\ref{prop:sparse}), the maximum achievable.
For any fixed coupling configuration between gyrostats and modes,
the number of subclasses grows as $2^{6K}$, where $K$ is the number
of gyrostats, making exhaustive characterization intractable for large
$K$; in contrast, non-canonical Hamiltonian structure (if and where
it arises from constraints on parameters) provides a systematic alternative
to finding invariants.

We identify conditions on the parameters, especially the nonlinear
coefficients, under which a GLOM admits a non-canonical Hamiltonian
structure, where its dynamics can be expressed in terms of a skew-symmetric
Poisson matrix $\mathrm{J}$ that satisfies the Jacobi identity. When
these conditions hold, the quadratic invariants beyond energy arise
as Casimir functions of $\mathrm{J}$ and can be read directly from
its nullspace. For physically motivated coupled hierarchies, including
models of 2D and 3D Rayleigh-B�nard convection, these conditions are
either satisfied exactly or otherwise impose explicit constraints
on the gyrostat parameters.

For nested Hamiltonian hierarchies obtained by progressively adding
gyrostats in a consistent pattern, the conditions for Hamiltonian
structure at each new level satisfy a recurrence (Proposition~\ref{prop:recurse}).
Furthermore, Casimir gradients project consistently under restriction
to lower-dimensional subspaces of the hierarchy (Theorem~\ref{prop:consist}),
so that invariants are compatible across models of increasing complexity.
This consistency property provides a geometric criterion for constructing
physically faithful reduced-order model hierarchies where the tools
of Hamiltonian dynamics can be deployed.

\section{Models and methods}

\subsection{GLOMs: setup and notation}

GLOMs are systems of $K$ gyrostats \citep{Gluhovsky1999}, where
each gyrostat $k\in\{1,\ldots,K\}$ involves three coupled modes 
\begin{align}
\dot{y}_{1}^{\left(k\right)} & =p_{k}y_{2}^{\left(k\right)}y_{3}^{\left(k\right)}+b_{k}y_{3}^{\left(k\right)}-c_{k}y_{2}^{\left(k\right)}\nonumber \\
\dot{y}_{2}^{\left(k\right)} & =q_{k}y_{3}^{\left(k\right)}y_{1}^{\left(k\right)}+c_{k}y_{1}^{\left(k\right)}-a_{k}y_{3}^{\left(k\right)}\nonumber \\
\dot{y}_{3}^{\left(k\right)} & =r_{k}y_{1}^{\left(k\right)}y_{2}^{\left(k\right)}+a_{k}y_{2}^{\left(k\right)}-b_{k}y_{1}^{\left(k\right)},\label{eq:1}
\end{align}
with real parameters $\{a_{k},b_{k},c_{k},p_{k},q_{k},r_{k}\}$ subject
to the energy conservation constraint $p_{k}+q_{k}+r_{k}=0$. The
$K$ gyrostats together couple $M\geq3$ modes, with the $k$th gyrostat
involving mode indices $m_{1}^{\left(k\right)},m_{2}^{\left(k\right)},m_{3}^{\left(k\right)}\in\{1,\ldots,M\}$.
The evolution of the $M$ mode amplitudes $x_{i}$ is given by superposition
of gyrostat contributions: $\dot{x}_{i}$ receives a contribution
$\dot{y}_{l}^{\left(k\right)}$ from Eq.~(\ref{eq:1}) whenever $m_{l}^{\left(k\right)}=i$,
with at most one equality possible for each $\left(i,k\right)$ pair.
Energy conservation $\dot{H}=0$ for $H=\frac{1}{2}\sum_{i=1}^{M}x_{i}^{2}$
follows directly from $p_{k}+q_{k}+r_{k}=0$ for each $k$, and holds
for all GLOMs regardless of any further parameter values.

We study the following five energy-conserving models.

$M=4$, $K=2$ (Model~1): 
\begin{align}
\dot{x}_{1} & =\left(p_{1}x_{2}x_{3}+b_{1}x_{3}-c_{1}x_{2}\right)\nonumber \\
\dot{x}_{2} & =\left(q_{1}x_{3}x_{1}+c_{1}x_{1}-a_{1}x_{3}\right)+\left\{ p_{2}x_{3}x_{4}+b_{2}x_{4}-c_{2}x_{3}\right\} \nonumber \\
\dot{x}_{3} & =\left(r_{1}x_{1}x_{2}+a_{1}x_{2}-b_{1}x_{1}\right)+\left\{ q_{2}x_{4}x_{2}+c_{2}x_{2}-a_{2}x_{4}\right\} \nonumber \\
\dot{x}_{4} & =\left\{ r_{2}x_{2}x_{3}+a_{2}x_{3}-b_{2}x_{2}\right\} ,\label{eq:5}
\end{align}

$M=5$, $K=2$ (Model~2): 
\begin{align}
\dot{x}_{1} & =\left(p_{1}x_{2}x_{3}+b_{1}x_{3}-c_{1}x_{2}\right)\nonumber \\
\dot{x}_{2} & =\left(q_{1}x_{3}x_{1}+c_{1}x_{1}-a_{1}x_{3}\right)\nonumber \\
\dot{x}_{3} & =\left(r_{1}x_{1}x_{2}+a_{1}x_{2}-b_{1}x_{1}\right)+\left\{ p_{2}x_{4}x_{5}+b_{2}x_{5}-c_{2}x_{4}\right\} \nonumber \\
\dot{x}_{4} & =\left\{ q_{2}x_{5}x_{3}+c_{2}x_{3}-a_{2}x_{5}\right\} \nonumber \\
\dot{x}_{5} & =\left\{ r_{2}x_{3}x_{4}+a_{2}x_{4}-b_{2}x_{3}\right\} ,\label{eq:6}
\end{align}

$M=5$, $K=3$ (Model~3): 
\begin{align}
\dot{x}_{1} & =\left(p_{1}x_{2}x_{3}+b_{1}x_{3}-c_{1}x_{2}\right)+\left[p_{3}x_{2}x_{4}+b_{3}x_{4}-c_{3}x_{2}\right]\nonumber \\
\dot{x}_{2} & =\left(q_{1}x_{3}x_{1}+c_{1}x_{1}-a_{1}x_{3}\right)+\left[q_{3}x_{4}x_{1}+c_{3}x_{1}-a_{3}x_{4}\right]\nonumber \\
\dot{x}_{3} & =\left(r_{1}x_{1}x_{2}+a_{1}x_{2}-b_{1}x_{1}\right)+\left\{ p_{2}x_{4}x_{5}+b_{2}x_{5}-c_{2}x_{4}\right\} \nonumber \\
\dot{x}_{4} & =\left\{ q_{2}x_{5}x_{3}+c_{2}x_{3}-a_{2}x_{5}\right\} +\left[r_{3}x_{1}x_{2}+a_{3}x_{2}-b_{3}x_{1}\right]\nonumber \\
\dot{x}_{5} & =\left\{ r_{2}x_{3}x_{4}+a_{2}x_{4}-b_{2}x_{3}\right\} .\label{eq:p7}
\end{align}

Different bracket styles identify contributions from gyrostat~1 (round),
gyrostat~2 (curly), and gyrostat~3 (square). Models~1 and~2 represent
the two natural ways to couple a second gyrostat to the first. In
Model~1 the two gyrostats share two modes ($x_{2}$ and $x_{3}$),
the densest possible two-gyrostat coupling, giving $M=K+2=4$. In
Model~2 they share a single mode ($x_{3}$), giving $M=2K+1=5$.
These configurations each generalise to arbitrary $K$: the dense
and sparse nested hierarchies have $M=K+2$ and $M=2K+1$ respectively,
as developed in Section~4. Model~3 illustrates a further effect:
a third gyrostat with modes $\{x_{1},x_{2},x_{4}\}$ is added to the
Model~2 configuration, creating a triangular cross-coupling in which
each gyrostat pair shares at~least one mode. Despite sharing $(K,M)=(3,5)$
with the $K=3$ instance of the dense hierarchy, Model~3 has a different
coupling structure and, as shown in Section~3, a lower maximum invariant
count.

\subsection{Finding quadratic invariants: the standard approach}

The search for invariants is generally a problem in linear algebra,
as illustrated by the single gyrostat with three modes 
\begin{align}
\dot{x}_{1} & =p_{1}x_{2}x_{3}+b_{1}x_{3}-c_{1}x_{2}\nonumber \\
\dot{x}_{2} & =q_{1}x_{3}x_{1}+c_{1}x_{1}-a_{1}x_{3}\nonumber \\
\dot{x}_{3} & =r_{1}x_{1}x_{2}+a_{1}x_{2}-b_{1}x_{1}\label{eq:3}
\end{align}
with $p_{1}+q_{1}+r_{1}=0$, for which we seek quadratic invariants
of the form $C_{1}=\sum_{i}\frac{1}{2}d_{i}x_{i}^{2}+\sum_{i<j}e_{ij}x_{i}x_{j}+\sum_{i}f_{i}x_{i}$
by solving $\dot{C}_{1}=0$. In this example, invariants acquire a
preferred coordinate system and $e_{ij}=0$ (Supplementary Information
(SI)), so that 
\begin{multline*}
\dot{C_{1}}=x_{1}x_{2}x_{3}\left(p_{1}d_{1}+q_{1}d_{2}+r_{1}d_{3}\right)+x_{1}x_{2}\left(-c_{1}d_{1}+c_{1}d_{2}+r_{1}f_{3}\right)\\
+x_{2}x_{3}\left(-a_{1}d_{2}+a_{1}d_{3}+p_{1}f_{1}\right)+x_{3}x_{1}\left(-b_{1}d_{3}+b_{1}d_{1}+q_{1}f_{2}\right)\\
+x_{1}\left(c_{1}f_{2}-b_{1}f_{3}\right)+x_{2}\left(a_{1}f_{3}-c_{1}f_{1}\right)+x_{3}\left(b_{1}f_{1}-a_{1}f_{2}\right)=0.
\end{multline*}
Setting the coefficient of each linearly independent term to zero
gives 
\begin{equation}
\left[\begin{array}{cccccc}
p_{1} & q_{1} & r_{1} & 0 & 0 & 0\\
-c_{1} & c_{1} & 0 & 0 & 0 & r_{1}\\
0 & -a_{1} & a_{1} & p_{1} & 0 & 0\\
b_{1} & 0 & -b_{1} & 0 & q_{1} & 0\\
0 & 0 & 0 & 0 & c_{1} & -b_{1}\\
0 & 0 & 0 & -c_{1} & 0 & a_{1}\\
0 & 0 & 0 & b_{1} & -a_{1} & 0
\end{array}\right]\left[\begin{array}{c}
d_{1}\\
d_{2}\\
d_{3}\\
f_{1}\\
f_{2}\\
f_{3}
\end{array}\right]=\left[\begin{array}{c}
0\\
0\\
0\\
0\\
0\\
0\\
0
\end{array}\right].\label{eq:p12}
\end{equation}
Denoting the above matrix as $\mathrm{A}$, the rank-nullity theorem
gives $\dim\left(\mathrm{NULL}\left(\mathrm{A}\right)\right)+\dim\left(\mathrm{RANGE}\left(\mathrm{A}\right)\right)=6$,
so the number of invariants is $6-\dim\left(\mathrm{RANGE}\left(\mathrm{A}\right)\right)$.
This is read from the column echelon form of $\mathrm{A}$, 
\[
\mathrm{A}_{ce}=\left[\begin{array}{cccccc}
p_{1} & 0 & 0 & 0 & 0 & 0\\
-c_{1} & -\frac{c_{1}r_{1}}{p_{1}} & 0 & 0 & 0 & 0\\
0 & -a_{1} & p_{1} & 0 & 0 & 0\\
b_{1} & -\frac{b_{1}q_{1}}{p_{1}} & 0 & -\frac{b_{1}q_{1}}{c_{1}} & 0 & 0\\
0 & 0 & 0 & -b_{1} & 0 & 0\\
0 & 0 & -c_{1} & 0 & 0 & 0\\
0 & 0 & b_{1} & \frac{a_{1}b_{1}}{c_{1}} & 0 & 0
\end{array}\right],
\]
which for the case with all parameters nonzero has four independent
columns, yielding two invariants (detailed calculations in SI). All
specialisations of the single gyrostat with some parameters zero also
have exactly two invariants \citep{Seshadri2023}, with the sequence
of column operations adjusted according to nonzero parameters. Since
$p_{1}+q_{1}+r_{1}=0$, energy $\frac{1}{2}\sum_{i=1}^{3}x_{i}^{2}$
is one invariant regardless of parameter values; the column echelon
form shows that a second invariant is a necessary consequence of this
energy conservation constraint (SI). Without this constraint, some
linear coefficients must vanish for any invariants at all to appear
\citep{Seshadri2023}.

\subsection{Non-canonical Hamiltonian structure of GLOMs}

\label{sec:ham}

The single gyrostat of Eq.~(\ref{eq:3}) admits a non-canonical Hamiltonian
formulation \citep{Gluhovsky2006} with Hamiltonian $H=\frac{1}{2}\sum_{i=1}^{3}x_{i}^{2}$
and skew-symmetric Poisson matrix 
\begin{equation}
\mathrm{J}=\left[\begin{array}{ccc}
0 & -c_{1} & p_{1}x_{2}+b_{1}\\
c_{1} & 0 & q_{1}x_{1}-a_{1}\\
-\left(p_{1}x_{2}+b_{1}\right) & -\left(q_{1}x_{1}-a_{1}\right) & 0
\end{array}\right],\label{eq:J_single}
\end{equation}
which recovers the vector field via $\dot{x}_{i}=\mathrm{J}_{ij}\frac{\partial H}{\partial x_{j}}=\mathrm{J}_{ij}x_{j}$
(repeated indices summed). For this to constitute a valid non-canonical
Hamiltonian system, $\mathrm{J}$ must satisfy the Jacobi condition
$\epsilon_{ijk}\mathrm{J}_{im}\frac{\partial\mathrm{J}_{jk}}{\partial x_{m}}=0$,
where $\epsilon_{ijk}$ is the alternating tensor. Unlike canonical
Hamiltonian mechanics \citep{Goldstein2002}, no constraint is placed
on the dimension of $\mathrm{J}$ or the structure of its entries
beyond skew-symmetry and, for example, odd-dimensional systems are
permitted. For the single gyrostat the Jacobi condition is identically
satisfied for all parameter values, making it unconditionally non-canonical
Hamiltonian \citep{Gluhovsky2006}.

When $\mathrm{J}$ is singular and a nontrivial nullspace vector is
the gradient of a scalar $C$, that scalar satisfies 
\begin{equation}
\mathrm{J}_{ij}\frac{\partial C}{\partial x_{j}}=0,\label{eq:casimir_cond}
\end{equation}
and is consequently conserved under the flow: 
\begin{equation}
\dot{C}=\frac{\partial C}{\partial x_{i}}\dot{x}_{i}=\frac{\partial C}{\partial x_{i}}\mathrm{J}_{ij}\frac{\partial H}{\partial x_{j}}=-\frac{\partial H}{\partial x_{i}}\mathrm{J}_{ij}\frac{\partial C}{\partial x_{j}}=0,\label{eq:casimir_cons}
\end{equation}
using skew-symmetry of $\mathrm{J}$ \citep{Shepherd1990}. Such scalars
$C$ are called Casimir functions. For the single gyrostat, $\mathrm{J}$
is skew-symmetric of odd order and therefore always singular; its
nullspace is spanned by $\left[a_{1}-q_{1}x_{1},\;b_{1}+p_{1}x_{2},\;c_{1}\right]^{T}$,
which is a gradient, giving the Casimir 
\[
C=-\tfrac{1}{2}q_{1}x_{1}^{2}+\tfrac{1}{2}p_{1}x_{2}^{2}+a_{1}x_{1}+b_{1}x_{2}+c_{1}x_{3}.
\]
This coincides with the second invariant found by the standard approach
and is the analogue of the squared angular momentum of the physical
gyroscope. In the geophysical context, Casimir functions of non-canonical
Hamiltonian fluid systems play a role analogous to integral invariants
such as enstrophy, potential vorticity, etc. \citep{Shepherd1990}:
they are conserved regardless of initial condition and confine trajectories
to submanifolds of the energy surface of the GLOM. Whether a GLOM
possesses Casimirs therefore has consequences that extend well beyond
the question of whether individual trajectories are regular or irregular,
as discussed in the Introduction.

For coupled GLOMs with $K>1$, Hamiltonian structure is not assured
and requires constraints on the model parameters. However, when the
Hamiltonian constraints hold, all quadratic invariants beyond energy
are recoverable directly from the nullspace of $\mathrm{J}$, whose
dimension is $M$, rather than from the nullspace of the larger matrix
$\mathrm{A}$ in Eq.~(\ref{eq:p12}). The Hamiltonian constraint
therefore both reduces the computational burden and provides a geometric
interpretation of the invariants. For any GLOM with $M$ modes coupled
by $K$ gyrostats, the skew-symmetric matrix satisfying $\dot{x}_{i}=\mathrm{J}_{ij}x_{j}$
is obtained as the superposition 
\begin{equation}
\mathrm{J}=\sum_{k=1}^{K}\mathrm{J}^{\left(k\right)},\label{eq:J_super}
\end{equation}
where $\mathrm{J}^{\left(k\right)}$ is the $M\times M$ matrix embedding
the $3\times3$ block 
\[
\mathrm{L}^{\left(k\right)}=\left[\begin{array}{ccc}
0 & -c_{k} & p_{k}x_{m_{2}^{\left(k\right)}}+b_{k}\\
c_{k} & 0 & q_{k}x_{m_{1}^{\left(k\right)}}-a_{k}\\
-\left(p_{k}x_{m_{2}^{\left(k\right)}}+b_{k}\right) & -\left(q_{k}x_{m_{1}^{\left(k\right)}}-a_{k}\right) & 0
\end{array}\right]
\]
at positions corresponding to the mode indices $m_{1}^{\left(k\right)},m_{2}^{\left(k\right)},m_{3}^{\left(k\right)}$
of the $k$th gyrostat, using the constraint $p_{k}+q_{k}+r_{k}=0$.
The Jacobi condition is then evaluated by substituting the elements
of $\mathrm{J}$ and collecting terms, and it imposes constraints
on the parameters that are present for each model in Section~4.

\smallskip{}
 \textit{Remark on the choice of\/ $\mathrm{J}$.} For each individual
gyrostat there are three valid choices of $\mathrm{L}^{\left(k\right)}$,
obtained by permuting which mode pair carries the nonlinear entry.
For $K$ gyrostats this yields $3^{K}$ candidate matrices $\mathrm{J}=\sum_{k}\mathrm{J}^{\left(k\right)}$.
The construction above selects one canonical assignment; others may
satisfy the Jacobi condition even when the canonical choice does not.
For the models 1-2 studied here we search for the Hamiltonian conditions
across all $3^{K}$ assignments (checked symbolically for $K\leq3$;
details in SI). \smallskip{}

\subsection{Nested and coupled hierarchies of Hamiltonian GLOMs}

\label{sec:hierarchies}

For GLOMs satisfying the Jacobi identity, we construct non-canonical
Hamiltonian hierarchies by progressively adding a gyrostat and deriving
the incremental conditions on parameters for the enlarged system to
remain Hamiltonian. We consider two classes of hierarchy.

\textit{Nested hierarchies.} Each increment of $K$ by one adds one
or two new modes. Two sub-cases arise. \textit{Sparse} nested hierarchies
extend the Model~2 pattern for $K=1,2,3,\ldots$ with $M=2K+1$:
at each step the new gyrostat shares one mode with its predecessor
and introduces two additional modes. \textit{Dense} nested hierarchies
extend the Model~1 pattern with $M=K+2$: the new gyrostat shares
two modes with the existing system and introduces one additional mode.
For both cases the incremental conditions on parameters for Hamiltonian
structure satisfy a recurrence (Proposition~\ref{prop:recurse}),
because the consistent coupling pattern propagates the Jacobi condition
locally between adjacent gyrostats.

\textit{Coupled hierarchies.} Here additional gyrostats couple existing
modes without introducing new ones; a simple recurrence does not hold
in general, but Hamiltonian conditions can still be derived. We investigate
two physically motivated examples from \citep{Gluhovsky2006}.

The first is the conservative core of a low-order model of 2D Rayleigh--B�nard
convection (Model~4). Its three gyrostats are all coupled through
a shared mode $x_{1}$, which represents the dominant convective cell,
giving $K=3$ and $M=7$: 
\begin{align}
\dot{x}_{1} & =-d_{1}x_{2}x_{3}-d_{2}x_{4}x_{5}-d_{3}x_{6}x_{7}\nonumber \\
\dot{x}_{2} & =d_{1}x_{3}x_{1}-d_{1}x_{3}\nonumber \\
\dot{x}_{3} & =d_{1}x_{2}\nonumber \\
\dot{x}_{4} & =d_{2}x_{5}x_{1}-d_{2}x_{5}\nonumber \\
\dot{x}_{5} & =d_{2}x_{4}\nonumber \\
\dot{x}_{6} & =d_{3}x_{7}x_{1}-d_{3}x_{7}\nonumber \\
\dot{x}_{7} & =d_{3}x_{6}.\label{eq:p38}
\end{align}

The second is the conservative core of a low-order model of 3D Rayleigh--B�nard
convection (Model~5), involving $K=5$ gyrostats and $M=8$ modes:
\begin{alignat}{1}
\dot{x}_{1} & =-x_{2}x_{3}-x_{4}x_{5}\nonumber \\
\dot{x}_{2} & =x_{3}x_{1}-x_{3}-\tfrac{1}{2}x_{5}x_{7}\nonumber \\
\dot{x}_{3} & =x_{2}\nonumber \\
\dot{x}_{4} & =x_{5}x_{1}-x_{5}-\tfrac{1}{2}x_{3}x_{7}\nonumber \\
\dot{x}_{5} & =x_{4}\nonumber \\
\dot{x}_{6} & =-2\beta x_{7}x_{8}\nonumber \\
\dot{x}_{7} & =2\beta x_{8}x_{6}-\beta x_{8}+\tfrac{1}{2}x_{3}x_{4}+\tfrac{1}{2}x_{5}x_{2}\nonumber \\
\dot{x}_{8} & =\beta x_{7}.\label{eq:p42}
\end{alignat}

For each hierarchy, symbolic computation is used to verify the Jacobi
condition at each level and to obtain gradient vectors of the Casimirs.
For nested hierarchies this is carried out for $K=1,2,3,4$; for coupled
hierarchies, among the $K!$ sub-hierarchies obtained by omitting
one gyrostat at a time, an illustrative ordering is analysed for each
model. All calculations are implemented in MATLAB using the Symbolic
Math Toolbox; code is provided in SI.

Figure 1 lists the hierarchy of GLOMs examined in this paper.

\begin{figure}
\includegraphics[scale=0.7]{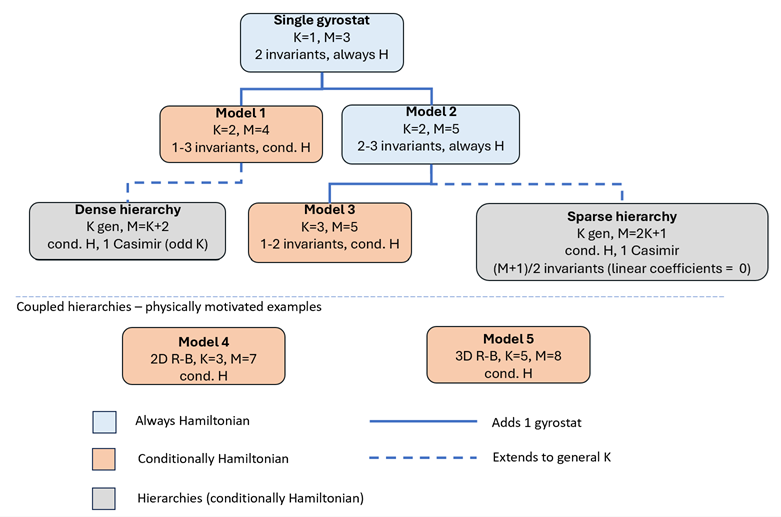}

\caption{Hierarchy of gyrostat low-order models (GLOMs) studied in this paper.
Each node shows the number of gyrostats $K$, mode count $M$, range
of quadratic invariant counts, and Hamiltonian status (H~=~Hamiltonian;
cond.~H~=~conditionally Hamiltonian, requiring parameter constraints).
Solid branches indicate that $K$ increases by one with a specific
coupling structure; dashed branches indicate that the $K=2$ model
is the base case of a general hierarchy examined for all $K$. Grey
nodes contain important illustrative (Dense and Sparse) hierarchies
of Section~4. Physically motivated coupled hierarchies modeling Rayleigh--B�nard
convection (Models~4--5) are also shown.}
\label{fig:hierarchy} 
\end{figure}

\section{Quadratic invariants by the standard approach}

\label{sec:standard}

\subsection{Illustrative calculations}

\subsubsection*{Model~1}

For the model of Eq.~(\ref{eq:5}), the mixed quadratic terms reduce
to a single nonzero coefficient $e_{14}$ (SI), so invariants take
the form 
\begin{multline}
C_{2}\!\left(x_{1},x_{2},x_{3},x_{4}\right)=\tfrac{1}{2}\!\left(d_{1}x_{1}^{2}+d_{2}x_{2}^{2}+d_{3}x_{3}^{2}+d_{4}x_{4}^{2}\right)+e_{14}x_{1}x_{4}\\
+f_{1}x_{1}+f_{2}x_{2}+f_{3}x_{3}+f_{4}x_{4}.\label{eq:inv_form1}
\end{multline}
Setting $\dot{C}_{2}=0$ and collecting linearly independent terms
yields the block system 
\[
\left[\begin{array}{cc}
\mathrm{A} & \mathrm{B}\\
\mathrm{0}_{4\times5} & \mathrm{D}
\end{array}\right]\left[\begin{array}{c}
\boldsymbol{\mu}\\
\boldsymbol{\nu}
\end{array}\right]=\mathbf{0},
\]
where $\boldsymbol{\mu}=\left[d_{1},d_{2},d_{3},d_{4},e_{14}\right]^{T}$,
$\boldsymbol{\nu}=\left[f_{1},f_{2},f_{3},f_{4}\right]^{T}$, and
\[
\mathrm{D}=\left[\begin{array}{cccc}
0 & c_{1} & -b_{1} & 0\\
-c_{1} & 0 & a_{1}+c_{2} & -b_{2}\\
b_{1} & -\left(a_{1}+c_{2}\right) & 0 & a_{2}\\
0 & b_{2} & -a_{2} & 0
\end{array}\right]
\]
is skew-symmetric. For the general case with all parameters nonzero,
$\mathrm{D}$ has full rank and therefore $\boldsymbol{\nu}=\mathbf{0}$;
the condition on invariants reduces to $\mathrm{A}\boldsymbol{\mu}=\mathbf{0}$,
with 
\begin{equation}
\mathrm{A}=\left[\begin{array}{ccccc}
p_{1} & q_{1} & r_{1} & 0 & r_{2}\\
-c_{1} & c_{1} & 0 & 0 & -b_{2}\\
0 & -\left(a_{1}+c_{2}\right) & \left(a_{1}+c_{2}\right) & 0 & 0\\
b_{1} & 0 & -b_{1} & 0 & a_{2}\\
0 & p_{2} & q_{2} & r_{2} & p_{1}\\
0 & 0 & -a_{2} & a_{2} & b_{1}\\
0 & b_{2} & 0 & -b_{2} & -c_{1}
\end{array}\right].\label{eq:A_model1}
\end{equation}
The column echelon form of $\mathrm{A}$ (SI) has four independent
columns, giving $\dim\!\left(\mathrm{NULL}\!\left(\mathrm{A}\right)\right)=1$
and a single invariant (energy) for the general case. In contrast,
the special case with all linear coefficients zero has three independent
invariants: 
\begin{multline}
C_{2}=\tfrac{1}{2}d_{1}\!\left(x_{1}-\tfrac{p_{1}}{r_{2}}x_{4}\right)^{\!2}+d_{2}\!\left(\tfrac{x_{2}^{2}}{2}-\tfrac{p_{2}-\frac{p_{1}q_{1}}{r_{2}}}{2r_{2}}x_{4}^{2}-\tfrac{q_{1}x_{1}x_{4}}{r_{2}}\right)\\
+d_{3}\!\left(\tfrac{x_{3}^{2}}{2}-\tfrac{q_{2}-\frac{p_{1}r_{1}}{r_{2}}}{2r_{2}}x_{4}^{2}-\tfrac{r_{1}x_{1}x_{4}}{r_{2}}\right),\label{eq:9}
\end{multline}
where the extra invariant arises from the linear dependence $\dot{x}_{1}=(p_{1}/r_{2})\dot{x}_{4}$
that holds when linear terms are absent. Energy is not independent
of the three quantities above. This model illustrates a general monotonicity
property:

\begin{proposition}[Monotonicity of invariant count]\label{prop:monotone}
Starting from a GLOM with some parameters set to zero, making any
subset of those parameters nonzero cannot increase the number of invariants.
\end{proposition} 
\begin{proof}
Let $\mathrm{A}_{0}$ denote the constraint matrix with certain parameters
set to zero, and let $s=\dim\!\left(\mathrm{RANGE}\!\left(\mathrm{A}_{0}\right)\right)$.
There exists an $s\times s$ submatrix $\mathrm{S}_{0}$ of $\mathrm{A}_{0}$
with nonzero determinant. Upon making the zero parameters nonzero
by amounts of order $\epsilon\ll1$, the matrix becomes $\mathrm{A}_{0}+\epsilon\mathrm{A}_{1}$,
and $\det\!\left(\mathrm{S}_{0}+\epsilon\mathrm{S}_{1}\right)=\det\mathrm{S}_{0}\!\left(1+\epsilon\,\mathrm{tr}\!\left(\mathrm{S}_{0}^{-1}\mathrm{S}_{1}\right)+O(\epsilon^{2})\right)\neq0$
for sufficiently small $\epsilon$. Hence $\dim\!\left(\mathrm{RANGE}\!\left(\mathrm{A}_{0}+\epsilon\mathrm{A}_{1}\right)\right)\geq s$,
and the null-space dimension cannot increase. 
\end{proof}
It follows that the minimum number of invariants for any GLOM configuration
is attained in the general case with all parameters nonzero, while
the maximum is attained in special cases with linear coefficients
absent. Table~\ref{tab:inv_counts} summarises the range for the
models studied.

\begin{table}[h]
\caption{Number of quadratic invariants for the models 1-3, with detailed calculations
in SI. The general case has all parameters nonzero; the no-feedback
case sets all linear coefficients $a_{k},b_{k},c_{k}$ to zero.}
\label{tab:inv_counts} \smallskip{}
 \centering %
\begin{tabular}{clcc}
\hline 
 & Model  & General case  & No linear feedback\tabularnewline
\hline 
1  & Single gyrostat, Eq.~(\ref{eq:3})  & 2  & 2\tabularnewline
2  & Model~1, Eq.~(\ref{eq:5})  & 1  & 3\tabularnewline
3  & Model~2, Eq.~(\ref{eq:6})  & 2  & 3\tabularnewline
4  & Model~3, Eq.~(\ref{eq:p7})  & 1  & 2\tabularnewline
\hline 
\end{tabular}
\end{table}

\subsubsection*{Model~2}

Invariants of this model take the form 
\[
C_{2}=\tfrac{1}{2}\sum_{i=1}^{5}d_{i}x_{i}^{2}+\sum_{i=1}^{5}f_{i}x_{i},
\]
since all mixed quadratic coefficients $e_{ij}$ vanish (SI). The
condition $\dot{C}_{2}=0$ yields the $13\times10$ constraint system
$\mathrm{A}\boldsymbol{u}=\mathbf{0}$ with $\boldsymbol{u}=\left[d_{1},\ldots,d_{5},f_{1},\ldots,f_{5}\right]^{T}$
and 
\[
\mathrm{A}=\left[\begin{array}{cccccccccc}
p_{1} & q_{1} & r_{1} & 0 & 0 & 0 & 0 & 0 & 0 & 0\\
-c_{1} & c_{1} & 0 & 0 & 0 & 0 & 0 & r_{1} & 0 & 0\\
0 & -a_{1} & a_{1} & 0 & 0 & p_{1} & 0 & 0 & 0 & 0\\
b_{1} & 0 & -b_{1} & 0 & 0 & 0 & q_{1} & 0 & 0 & 0\\
0 & 0 & p_{2} & q_{2} & r_{2} & 0 & 0 & 0 & 0 & 0\\
0 & 0 & -c_{2} & c_{2} & 0 & 0 & 0 & 0 & 0 & r_{2}\\
0 & 0 & 0 & -a_{2} & a_{2} & 0 & 0 & p_{2} & 0 & 0\\
0 & 0 & b_{2} & 0 & -b_{2} & 0 & 0 & 0 & q_{2} & 0\\
0 & 0 & 0 & 0 & 0 & 0 & c_{1} & -b_{1} & 0 & 0\\
0 & 0 & 0 & 0 & 0 & -c_{1} & 0 & a_{1} & 0 & 0\\
0 & 0 & 0 & 0 & 0 & b_{1} & -a_{1} & 0 & c_{2} & -b_{2}\\
0 & 0 & 0 & 0 & 0 & 0 & 0 & -c_{2} & 0 & a_{2}\\
0 & 0 & 0 & 0 & 0 & 0 & 0 & b_{2} & -a_{2} & 0
\end{array}\right].
\]
For the general case with all parameters nonzero, column reduction
(SI) gives $\dim\!\left(\mathrm{RANGE}\!\left(\mathrm{A}\right)\right)=8$,
so there are two quadratic invariants. In contrast, when all linear
coefficients are zero, the condition reduces to 
\begin{align}
d_{1}p_{1}+d_{2}q_{1}+d_{3}r_{1} & =0\nonumber \\
d_{3}p_{2}+d_{4}q_{2}+d_{5}r_{2} & =0,\label{eq:sparse_cond2}
\end{align}
two equations in five unknowns, giving three invariants.

\subsubsection*{Model~3}

Invariants again take the form $\tfrac{1}{2}\sum_{i}d_{i}x_{i}^{2}+\sum_{i}f_{i}x_{i}$
(mixed quadratic terms vanish, SI), with the number of invariants
determined by the nontrivial solutions to $\mathrm{A}\boldsymbol{u}=\mathbf{0}$,
where $\mathrm{A}\in\mathbb{R}^{16\times10}$ is given in SI. For
the general case without parameter restrictions, the null-space is
one-dimensional and only energy is invariant. In the special case
without linear feedbacks, the condition reduces to three equations,
\begin{align}
d_{1}p_{1}+d_{2}q_{1}+d_{3}r_{1} & =0\nonumber \\
d_{3}p_{2}+d_{4}q_{2}+d_{5}r_{2} & =0\nonumber \\
d_{1}p_{3}+d_{2}q_{3}+d_{4}r_{3} & =0,\label{eq:m3_nolinear}
\end{align}
in five unknowns, yielding two invariants. Comparing with Model~2
(Eq.~\ref{eq:sparse_cond2}), the addition of a third gyrostat whose
modes cross-couple the first two reduces the maximum invariant count
from three to two. This confirms that, with $M$ fixed, increasing
$K$ does not raise, and typically lowers, the invariant count.

\subsection{Sparse models without linear feedback}

\label{sec:sparse}

Consider the sparse hierarchy without any linear feedback terms present.
The pattern in Eq.~(\ref{eq:sparse_cond2} generalises to this hierarchy,
where each member of $K$ gyrostats has $M=2K+1$ modes: such, for
any $K$, the condition $\dot{C}=0$ reduces to $K$ equations $p_{k}d_{m_{1}^{(k)}}+q_{k}d_{m_{2}^{(k)}}+r_{k}d_{m_{3}^{(k)}}=0$
in $2K+1$ unknowns $d_{i}$ (linear and mixed quadratic terms vanish;
see Appendix~1 and SI). This immediately gives:

\begin{theorem}[Invariant count for sparse hierarchies]\label{prop:sparse}
Sparse GLOMs without linear feedback and with $M=2K+1$ modes have
exactly $\left(M+1\right)/2=K+1$ independent quadratic invariants.
\end{theorem} 
\begin{proof}
The $K$ constraint equations $p_{k}d_{m_{1}^{(k)}}+q_{k}d_{m_{2}^{(k)}}+r_{k}d_{m_{3}^{(k)}}=0$,
$k=1,\ldots,K$, have coefficient matrix of full rank for all $K\geq1$.
In the sparse hierarchy, the $k$-th gyrostat introduces the even-indexed
mode $x_{2k}$ for the first time: no other gyrostat involves $x_{2k}$,
so the variable $d_{2k}$ appears in the $k$-th constraint equation
and in no other. Ordering the $K$ equations by $k=1,\ldots,K$ and
the $2K+1$ unknowns so that $d_{2},d_{4},\ldots,d_{2K}$ come first,
the $K\times K$ submatrix formed by selecting only the $K$ columns
$\{d_{2},d_{4},\ldots,d_{2K}\}$ is diagonal, with diagonal entries
$q_{k}$ (for the canonical representation in which $x_{2k}$ occupies
the second position of the $k$-th gyrostat triple). Since $q_{k}\neq0$
for generic parameters, this $K\times K$ submatrix has rank $K$,
so the full $K\times(2K+1)$ constraint matrix also has rank $K$,
giving null-space dimension $M-K=\left(2K+1\right)-K=K+1=\left(M+1\right)/2$. 
\end{proof}
The explicit invariants are listed in Table~\ref{tab:sparse_inv}
for $K=1,2,3,4$. Each $C_{K,m}$ (other than the first and the last)
involves only the squared amplitudes of the odd-indexed modes $x_{1},x_{3},\ldots,x_{2K-1}$
together with one even-indexed mode $x_{2m}$, with coefficients determined
by ratios of the nonlinear parameters (recall that the linear parameters
are zero). For each $K$, the invariants sum to $\sum_{i=1}^{M}x_{i}^{2}$,
and adding a gyrostat (incrementing $K$ by one, introducing modes
$x_{2K}$ and $x_{2K+1}$) preserves all but the last invariant of
the previous level in the hierarchy.

\begin{table}[h]
\caption{Quadratic invariants for sparse GLOMs without any linear feedback
terms being present ($M=2K+1$). Calculation and verification for
each case are in SI.}
\label{tab:sparse_inv} \smallskip{}
 \centering %
\begin{tabular}{ccp{9.5cm}}
\hline 
$K$  & $M$  & Invariants $C_{K,m}$, $m=1,\ldots,K+1$\tabularnewline
\hline 
$1$  & $3$  & $-\dfrac{q_{1}}{p_{1}}x_{1}^{2}+x_{2}^{2}$,\quad{}$-\dfrac{r_{1}}{p_{1}}x_{1}^{2}+x_{3}^{2}$\tabularnewline
$2$  & $5$  & $-\dfrac{q_{1}}{p_{1}}x_{1}^{2}+x_{2}^{2}$,\quad{}$\dfrac{q_{1}r_{1}}{p_{1}p_{2}}x_{1}^{2}-\dfrac{q_{2}}{p_{2}}x_{3}^{2}+x_{4}^{2}$,\quad{}$\dfrac{r_{1}r_{2}}{p_{1}p_{2}}x_{1}^{2}-\dfrac{r_{2}}{p_{2}}x_{3}^{2}+x_{5}^{2}$\tabularnewline
$3$  & $7$  & $-\dfrac{q_{1}}{p_{1}}x_{1}^{2}+x_{2}^{2}$,\quad{}$\dfrac{q_{2}r_{1}}{p_{1}p_{2}}x_{1}^{2}-\dfrac{q_{2}}{p_{2}}x_{3}^{2}+x_{4}^{2}$,\quad{}$-\dfrac{q_{3}r_{1}r_{2}}{p_{1}p_{2}p_{3}}x_{1}^{2}+\dfrac{q_{3}r_{2}}{p_{2}p_{3}}x_{3}^{2}-\dfrac{q_{3}}{p_{3}}x_{5}^{2}+x_{6}^{2}$,\quad{}$-\dfrac{r_{1}r_{2}r_{3}}{p_{1}p_{2}p_{3}}x_{1}^{2}+\dfrac{r_{2}r_{3}}{p_{2}p_{3}}x_{3}^{2}-\dfrac{r_{3}}{p_{3}}x_{5}^{2}+x_{7}^{2}$\tabularnewline
$4$  & $9$  & $-\dfrac{q_{1}}{p_{1}}x_{1}^{2}+x_{2}^{2}$,\quad{}$\dfrac{q_{2}r_{1}}{p_{1}p_{2}}x_{1}^{2}-\dfrac{q_{2}}{p_{2}}x_{3}^{2}+x_{4}^{2}$,\quad{}$-\dfrac{q_{3}r_{1}r_{2}}{p_{1}p_{2}p_{3}}x_{1}^{2}+\dfrac{q_{3}r_{2}}{p_{2}p_{3}}x_{3}^{2}-\dfrac{q_{3}}{p_{3}}x_{5}^{2}+x_{6}^{2}$,\quad{}$\dfrac{q_{4}r_{1}r_{2}r_{3}}{p_{1}p_{2}p_{3}p_{4}}x_{1}^{2}-\dfrac{q_{4}r_{2}r_{3}}{p_{2}p_{3}p_{4}}x_{3}^{2}+\dfrac{q_{4}r_{3}}{p_{3}p_{4}}x_{5}^{2}-\dfrac{q_{4}}{p_{4}}x_{7}^{2}+x_{8}^{2}$,\quad{}$\dfrac{r_{1}r_{2}r_{3}r_{4}}{p_{1}p_{2}p_{3}p_{4}}x_{1}^{2}-\dfrac{r_{2}r_{3}r_{4}}{p_{2}p_{3}p_{4}}x_{3}^{2}+\dfrac{r_{3}r_{4}}{p_{3}p_{4}}x_{5}^{2}-\dfrac{r_{4}}{p_{4}}x_{7}^{2}+x_{9}^{2}$\tabularnewline
\hline 
\end{tabular}
\end{table}

\subsection{Sensitivity of invariant count to model configuration}

\label{sec:sensitivity}

The sparse models in the special case without any linear feedback
terms are exceptional in admitting a straightforward characterisation
of invariants across the hierarchy. For general GLOMs, even fixing
$(K,M)$ as well as the structural coupling between gyrostats, the
number of invariants depends sensitively on which linear feedbacks
are present. We illustrate this for Models~1 and~2, then identify
the factors behind the difficulty of extending the approach to larger
models.

\subsubsection*{Model~1}

With all nonlinear coefficients $p_{i},q_{i},r_{i}$ nonzero, the
invariant count is controlled entirely by the four linear feedback
parameters $\{b_{1},c_{1},a_{2},b_{2}\}$. Table~\ref{tab:model1_cond}
summarises the conditions. Three invariants occur only when all four
feedbacks are simultaneously zero, in which case the vector field
satisfies $r_{2}\dot{x}_{1}-p_{1}\dot{x}_{4}=0$, making two components
linearly dependent and yielding an invariant of the form $\left(x_{1}-p_{1}x_{4}/r_{2}\right)^{2}$.
Two invariants arise whenever the nonzero feedbacks avoid all four
forbidden pairs $(b_{1},c_{1})$, $(b_{1},b_{2})$, $(c_{1},a_{2})$,
$(a_{2},b_{2})$. All remaining cases, including the general case
with all feedbacks nonzero, yield only energy as the single invariant.

\begin{table}[h]
\caption{Conditions on linear feedback parameters for each invariant count
in Model~1, assuming all nonlinear coefficients $p_{i},q_{i},r_{i}$
are nonzero.}
\label{tab:model1_cond} \smallskip{}
 \centering %
\begin{tabular}{cp{10.5cm}}
\hline 
Invariants  & Conditions on $\{b_{1},c_{1},a_{2},b_{2}\}$\tabularnewline
\hline 
$3$  & $b_{1}=c_{1}=a_{2}=b_{2}=0$. Results in linear dependence $r_{2}\dot{x}_{1}=p_{1}\dot{x}_{4}$;
energy is not independent of the three invariants.\tabularnewline
$2$  & No forbidden pair is simultaneously nonzero, i.e.\ not $\left(b_{1}\neq0\text{ and }c_{1}\neq0\right)$,
not $\left(b_{1}\neq0\text{ and }b_{2}\neq0\right)$, not $\left(c_{1}\neq0\text{ and }a_{2}\neq0\right)$,
not $\left(a_{2}\neq0\text{ and }b_{2}\neq0\right)$. There are exactly
six non-empty subsets satisfying this: $\{b_{1}\}$, $\{b_{1},a_{2}\}$,
$\{a_{2}\}$, $\{b_{2}\}$, $\{c_{1}\}$, $\{c_{1},b_{2}\}$.\tabularnewline
$1$  & All other cases; in particular the general case with all feedbacks
nonzero.\tabularnewline
\hline 
\end{tabular}
\end{table}

The sensitivity to fine structure of the feedback terms is illustrated
by contrasting two subclasses that each have four nonzero linear parameters
(two from each gyrostat). The system 
\begin{align}
\dot{x}_{1} & =\left(p_{1}x_{2}x_{3}-c_{1}x_{2}\right)\nonumber \\
\dot{x}_{2} & =\left(q_{1}x_{3}x_{1}+c_{1}x_{1}-a_{1}x_{3}\right)+\left\{ p_{2}x_{3}x_{4}+b_{2}x_{4}-c_{2}x_{3}\right\} \nonumber \\
\dot{x}_{3} & =\left(r_{1}x_{1}x_{2}+a_{1}x_{2}\right)+\left\{ q_{2}x_{4}x_{2}+c_{2}x_{2}\right\} \nonumber \\
\dot{x}_{4} & =\left\{ r_{2}x_{2}x_{3}-b_{2}x_{2}\right\} \label{eq:p20}
\end{align}
has two invariants (nonzero $c_{1}$ and $b_{2}$ avoid all forbidden
pairs), while 
\begin{align}
\dot{x}_{1} & =\left(p_{1}x_{2}x_{3}-c_{1}x_{2}\right)\nonumber \\
\dot{x}_{2} & =\left(q_{1}x_{3}x_{1}+c_{1}x_{1}-a_{1}x_{3}\right)+\left\{ p_{2}x_{3}x_{4}-c_{2}x_{3}\right\} \nonumber \\
\dot{x}_{3} & =\left(r_{1}x_{1}x_{2}+a_{1}x_{2}\right)+\left\{ q_{2}x_{4}x_{2}+c_{2}x_{2}-a_{2}x_{4}\right\} \nonumber \\
\dot{x}_{4} & =\left\{ r_{2}x_{2}x_{3}+a_{2}x_{3}\right\} \label{eq:p21}
\end{align}
has only energy as invariant (nonzero $c_{1}$ and $a_{2}$ form the
forbidden pair $(c_{1},a_{2})$). The second model has effectively
three-dimensional conservative dynamics on the energy surface, making
irregular trajectories possible without any forcing or dissipation.
Neither model has symmetries that could predict the difference between
them in this respect and, moreover, a mixed quadratic term $x_{1}x_{4}$
is present in every non-energy invariant of Model~1 (SI).

\subsubsection*{Model~2}

If all nonlinear coefficients are nonzero, the invariant count depends
only on two linear feedback parameters: three invariants arise if
and only if $c_{1}=a_{2}=0$, giving the system 
\begin{align*}
\dot{x}_{1} & =\left(p_{1}x_{2}x_{3}+b_{1}x_{3}\right)\\
\dot{x}_{2} & =\left(q_{1}x_{3}x_{1}-a_{1}x_{3}\right)\\
\dot{x}_{3} & =\left(r_{1}x_{1}x_{2}+a_{1}x_{2}-b_{1}x_{1}\right)+\left\{ p_{2}x_{4}x_{5}+b_{2}x_{5}-c_{2}x_{4}\right\} \\
\dot{x}_{4} & =\left\{ q_{2}x_{5}x_{3}+c_{2}x_{3}\right\} \\
\dot{x}_{5} & =\left\{ r_{2}x_{3}x_{4}-b_{2}x_{3}\right\} .
\end{align*}
All other subclasses (with $c_{1}\neq0$ or $a_{2}\neq0$ or both)
retain exactly two invariants. Two invariants in Model 2 still permit
irregular dynamics but three invariants would not. Comparing with
Model~1 (Table~\ref{tab:model1_cond}), it is seen that the sparser
coupling of Model~2 is more robust since adding the second gyrostat
by introducing two new modes preserves the invariant count of the
single gyrostat in general, without restrictions on parameters, while
in Model~1 it reduces the number of invariants from two to one unless
parameters are restricted.

\subsubsection*{Limitations of the standard approach}

The preceding calculations reveal various reasons why the standard
approach cannot be extended systematically to larger GLOMs.

First, constructing the constraint matrix becomes rapidly challenging
as $\mathrm{A}$ grows rapidly with $M$. Identifying which mixed
quadratic terms vanish for all subclasses of a given configuration
requires specialized analysis for each model; without such reductions
the number of columns of $\mathrm{A}$ grows quadratically with $M$,
and the column reduction becomes unwieldy. Second, the number of parameter
subclasses grows as $2^{6K}$, making exhaustive characterisation
exponentially expensive as $K$ increases; the sensitive dependence
of invariant count on the exact configuration (as illustrated above)
means there does not appear to be an evident shortcut through this
combinatorially large space. Third, degenerate cases, where components
of the vector field are linearly dependent, can cause the null-space
dimension to overcount the true invariant count. For example, the
subclass of Eq.~(\ref{eq:5}) with $p_{1}=b_{1}=c_{1}=0$ satisfies
$\dot{x}_{1}=0$; the null-space of $\mathrm{A}$ indicates four invariants,
but $x_{1}$ and $x_{1}^{2}$ are not independent, giving only three
independent invariants.

\subsection{Limited role for symmetries}

Symmetries of the vector field can simplify the analysis of invariants
(Appendix~2), but are of limited help for the general problem. The
presence of symmetries in GLOMs requires many linear coefficients
to vanish simultaneously, so they arise only in special subclasses.
Moreover, symmetries are sufficient but not necessary for coefficients
to vanish: when $e_{ij}=0$ in the general case, this follows from
the independence of the constraint equations rather than any underlying
symmetry.

To take an example, the Euler gyrostat 
\begin{align*}
\dot{x}_{1} & =p_{1}x_{2}x_{3}\\
\dot{x}_{2} & =q_{1}x_{3}x_{1}\\
\dot{x}_{3} & =r_{1}x_{1}x_{2}
\end{align*}
is invariant under the three sign-flip transformations $\left(x_{1},x_{2},x_{3}\right)\to\left(-x_{1},-x_{2},x_{3}\right)$,
$\left(x_{1},-x_{2},-x_{3}\right)$, $\left(-x_{1},x_{2},-x_{3}\right)$.
Any invariant must be maintained under all three, which immediately
rules out mixed quadratic terms $x_{i}x_{j}$ and linear terms $x_{i}$,
leaving only diagonal quadratic invariants (Table~\ref{tab:inv_counts}).
The sparse model in Eq. (\ref{eq:6}) without any linear feedback
terms has seven such sign-flip symmetries (SI), which similarly restrict
invariants to diagonal quadratic form, consistent with Theorem~\ref{prop:sparse}.
Sparse hierarchies without linear feedback share this structure because
their symmetry group is generated by compositions of the single Euler-gyrostat
symmetries. In contrast, models with linear feedback terms break these
symmetries, and the invariant count must then be determined by the
full constraint analysis of the above sections.

The above analysis shows that while the standard approach succeeds
for individual cases, it cannot be systematically extended to the
full class of GLOMs. We turn in the next section to non-canonical
Hamiltonian structure, which provides a framework that scales more
readily with the hierarchy.

\section{Invariants from Hamiltonian structure}

\label{sec:results_ham}

The standard approach of Section~\ref{sec:standard} succeeds for
particular models of small or moderate size but does not scale. We
now exploit the non-canonical Hamiltonian framework introduced in
Section~\ref{sec:ham}: for each model we evaluate the Jacobi condition
on $\mathrm{J}=\sum_{k}\mathrm{J}^{(k)}$, identify the parameter
constraints that make it zero, and identify the Casimir functions
from the resulting nullspace of $\mathrm{J}$. Two recurring themes
emerge. First, Hamiltonian structure typically requires constraints
on the nonlinear coefficients, not only the linear feedbacks, because
only the nonlinear entries of $\mathrm{J}^{(k)}$ produce nonzero
derivatives $\partial\mathrm{J}_{jk}/\partial x_{m}$. Second, restricting
to Hamiltonian subclasses markedly simplifies the characterisation
of quadratic invariants, as Table~\ref{tab:model1_ham} illustrates.

\subsection{Hamiltonian conditions and Casimirs: Models~1--3}

\subsubsection*{Model~1}

Substituting the superposition $\mathrm{J}=\mathrm{J}^{(1)}+\mathrm{J}^{(2)}$
into Eq.~(\ref{eq:5}), 
\begin{multline}
\mathrm{J}=\left[\begin{array}{cccc}
0 & -c_{1} & p_{1}x_{2}+b_{1} & 0\\
c_{1} & 0 & q_{1}x_{1}-a_{1} & 0\\
-\left(p_{1}x_{2}+b_{1}\right) & -\left(q_{1}x_{1}-a_{1}\right) & 0 & 0\\
0 & 0 & 0 & 0
\end{array}\right]+\left[\begin{array}{cccc}
0 & 0 & 0 & 0\\
0 & 0 & -c_{2} & p_{2}x_{3}+b_{2}\\
0 & c_{2} & 0 & q_{2}x_{2}-a_{2}\\
0 & -\left(p_{2}x_{3}+b_{2}\right) & -\left(q_{2}x_{2}-a_{2}\right) & 0
\end{array}\right]\\
=\left[\begin{array}{cccc}
0 & -c_{1} & p_{1}x_{2}+b_{1} & 0\\
c_{1} & 0 & q_{1}x_{1}-a_{1}-c_{2} & p_{2}x_{3}+b_{2}\\
-\left(p_{1}x_{2}+b_{1}\right) & -\left(q_{1}x_{1}-a_{1}\right)+c_{2} & 0 & q_{2}x_{2}-a_{2}\\
0 & -\left(p_{2}x_{3}+b_{2}\right) & -\left(q_{2}x_{2}-a_{2}\right) & 0
\end{array}\right].\label{eq:p25}
\end{multline}
Expanding terms with nonzero $\partial\mathrm{J}_{jk}/\partial x_{m}$,
the Jacobi condition reduces to 
\begin{equation}
\epsilon_{ijk}\mathrm{J}_{im}\frac{\partial\mathrm{J}_{jk}}{\partial x_{m}}=2p_{1}p_{2}\left(x_{2}-x_{3}\right)+2\left(p_{2}b_{1}-p_{1}b_{2}-q_{2}c_{1}\right)=0.\label{eq:jac1}
\end{equation}
This is satisfied for all $\mathbf{x}$ if and only if both 
\begin{equation}
p_{1}p_{2}=0\qquad\text{and}\qquad p_{2}b_{1}-p_{1}b_{2}-q_{2}c_{1}=0.\label{eq:ham1_cond}
\end{equation}
Since both $p_{1}=0$ and $p_{2}=0$ simultaneously would leave neither
gyrostat with three nonlinear terms, we restrict our treatment to
cases where exactly one of them vanishes.

\textit{Case~1: $p_{1}=0$.} The Hamiltonian constraint then also
requires $b_{1}=c_{1}=0$, yielding $\dot{x}_{1}=0$ identically.
This is a degenerate model: $x_{1}$ is a constant of motion, and
the system reduces to a single gyrostat in the variables $x_{2},x_{3},x_{4}$.
Both columns of $\mathrm{NULL}\left(\mathrm{J}\right)$ are gradient
vectors, giving Casimirs $C_{a}=x_{1}$ and 
\[
C_{a}=-\tfrac{1}{2}q_{2}x_{2}^{2}+\tfrac{1}{2}p_{2}x_{3}^{2}-q_{1}x_{1}x_{4}+a_{2}x_{2}+b_{2}x_{3}+\left(a_{1}+c_{2}\right)x_{4}.
\]
Note that $\mathrm{NULL}\left(\mathrm{A}\right)$ additionally contains
$x_{1}^{2}$, which is not independent of the two Casimirs above.

\textit{Case~2: $p_{2}=0$, additionally $c_{1}=b_{2}=0$.} This
nondegenerate case gives the GLOM 
\begin{align*}
\dot{x}_{1} & =\left(p_{1}x_{2}x_{3}+b_{1}x_{3}\right)\\
\dot{x}_{2} & =\left(q_{1}x_{3}x_{1}-a_{1}x_{3}\right)+\left\{ -c_{2}x_{3}\right\} \\
\dot{x}_{3} & =\left(r_{1}x_{1}x_{2}+a_{1}x_{2}-b_{1}x_{1}\right)+\left\{ q_{2}x_{4}x_{2}+c_{2}x_{2}-a_{2}x_{4}\right\} \\
\dot{x}_{4} & =\left\{ r_{2}x_{2}x_{3}+a_{2}x_{3}\right\} .
\end{align*}
The corresponding $\mathrm{J}$ has $\det\mathrm{J}=0$ and $\mathrm{rank}\left(\mathrm{J}\right)=2$.
The nullspace $\mathrm{NULL}\left(\mathrm{J}\right)$ is spanned by
two vectors, of which only 
\[
\mathbf{v}_{1}=\left[\begin{array}{c}
a_{1}+c_{2}-q_{1}x_{1}\\
b_{1}+p_{1}x_{2}\\
0\\
0
\end{array}\right]
\]
is a gradient vector, yielding the Casimir 
\begin{equation}
C_{a}=-\frac{q_{1}}{2}x_{1}^{2}+\frac{p_{1}}{2}x_{2}^{2}+\left(a_{1}+c_{2}\right)x_{1}+b_{1}x_{2}.\label{eq:cas_m1}
\end{equation}
The second basis vector of $\mathrm{NULL}\left(\mathrm{J}\right)$
is not a gradient, so there is only one Casimir.

Table~\ref{tab:model1_ham} summarises the invariant-count conditions
for Case~2, replacing the detailed subclass enumeration of the non-Hamiltonian
case. The key contrast with the general (non-Hamiltonian) analysis
of Table~\ref{tab:model1_cond} is that within the Hamiltonian subclass
the conditions are governed by only two parameters, $b_{1}$ and $a_{2}$,
rather than all four linear feedbacks.

\begin{table}[h]
\caption{Invariant-count conditions for the nondegenerate Hamiltonian subclass
of Model~1 ($p_{2}=c_{1}=b_{2}=0$, all remaining nonlinear coefficients
nonzero).}
\label{tab:model1_ham} \smallskip{}
 \centering %
\begin{tabular}{cp{10.5cm}}
\hline 
Invariants  & Conditions\tabularnewline
\hline 
$2$  & General Hamiltonian case: $b_{1}\neq0$ or $a_{2}\neq0$. One Casimir
$C_{a}$ given by Eq.~(\ref{eq:cas_m1}); energy $H$ is the second
invariant.\tabularnewline
$3$  & $b_{1}=a_{2}=0$: the vector field satisfies $r_{2}\dot{x}_{1}=p_{1}\dot{x}_{4}$,
creating a linear dependence and an additional invariant $(x_{1}-p_{1}x_{4}/r_{2})^{2}$.
In this degenerate case energy is not independent of the three invariants.\tabularnewline
\hline 
\end{tabular}
\end{table}

\subsubsection*{Model~2}

The Poisson matrix for Eq.~(\ref{eq:6}) is 
\[
\mathrm{J}=\left[\begin{array}{ccccc}
0 & -c_{1} & p_{1}x_{2}+b_{1} & 0 & 0\\
c_{1} & 0 & q_{1}x_{1}-a_{1} & 0 & 0\\
-\left(p_{1}x_{2}+b_{1}\right) & -\left(q_{1}x_{1}-a_{1}\right) & 0 & -c_{2} & p_{2}x_{4}+b_{2}\\
0 & 0 & c_{2} & 0 & q_{2}x_{3}-a_{2}\\
0 & 0 & -\left(p_{2}x_{4}+b_{2}\right) & -\left(q_{2}x_{3}-a_{2}\right) & 0
\end{array}\right],
\]
and the Jacobi condition simplifies to 
\begin{equation}
\epsilon_{ijk}\mathrm{J}_{im}\frac{\partial\mathrm{J}_{jk}}{\partial x_{m}}=2q_{2}\left(q_{1}x_{1}+p_{1}x_{2}+b_{1}-a_{1}\right)=0.\label{eq:jac2}
\end{equation}
Since the factor $q_{1}x_{1}+p_{1}x_{2}+b_{1}-a_{1}$ cannot vanish
identically unless $p_{1}=q_{1}=b_{1}=a_{1}=0$ (excluded by the requirement
of at least five nonzero nonlinear coefficients), the sole Hamiltonian
model we examine has 
\begin{equation}
q_{2}=0.\label{eq:ham2_cond}
\end{equation}
With $q_{2}=0$, the model becomes 
\begin{align}
\dot{x}_{1} & =\left(p_{1}x_{2}x_{3}+b_{1}x_{3}-c_{1}x_{2}\right)\nonumber \\
\dot{x}_{2} & =\left(q_{1}x_{3}x_{1}+c_{1}x_{1}-a_{1}x_{3}\right)\nonumber \\
\dot{x}_{3} & =\left(r_{1}x_{1}x_{2}+a_{1}x_{2}-b_{1}x_{1}\right)+\left\{ p_{2}x_{4}x_{5}+b_{2}x_{5}-c_{2}x_{4}\right\} \nonumber \\
\dot{x}_{4} & =\left\{ c_{2}x_{3}-a_{2}x_{5}\right\} \nonumber \\
\dot{x}_{5} & =\left\{ -p_{2}x_{3}x_{4}+a_{2}x_{4}-b_{2}x_{3}\right\} .\label{eq:m2_ham}
\end{align}
For all other parameters nonzero, $\mathrm{NULL}\left(\mathrm{J}\right)$
is spanned by the single gradient vector 
\begin{equation}
\left[\begin{array}{c}
a_{2}\left(a_{1}-q_{1}x_{1}\right)\\
a_{2}\left(b_{1}+p_{1}x_{2}\right)\\
a_{2}c_{1}\\
c_{1}\left(b_{2}+p_{2}x_{4}\right)\\
c_{1}c_{2}
\end{array}\right],\label{eq:cas_m2_grad}
\end{equation}
yielding the Casimir 
\begin{equation}
C_{a}=-\tfrac{1}{2}a_{2}q_{1}x_{1}^{2}+\tfrac{1}{2}a_{2}p_{1}x_{2}^{2}+\tfrac{1}{2}c_{1}p_{2}x_{4}^{2}+a_{1}a_{2}x_{1}+a_{2}b_{1}x_{2}+a_{2}c_{1}x_{3}+c_{1}b_{2}x_{4}+c_{1}c_{2}x_{5}.\label{eq:cas_m2}
\end{equation}
This is confirmed by the null-space of $\mathrm{A}$ (SI). The above
model, and associated Poisson matrix $\mathrm{J}$, can be further
specialized by restricting parameters. For example $c_{1}=a_{2}=0$
yields two independent Casimirs. This restriction also increases invariant
count for the non-Hamiltonian cases, so such reductions are not unique
to the non-Hamiltonian cases. The difference is, firstly, the avoidance
of spurious degeneracies in the Hamiltonian models and, secondly,
the legibility of the invariants derived within the Hamiltonian framework
and its role in generalizing to higher $K$ and $M$.

\subsubsection*{Model~3}

With the addition of a third gyrostat the Poisson matrix becomes 
\[
\mathrm{J}=\left[\begin{array}{ccccc}
0 & -c_{1}-c_{3} & p_{1}x_{2}+b_{1} & b_{3}+p_{3}x_{2} & 0\\
c_{1}+c_{3} & 0 & q_{1}x_{1}-a_{1} & q_{3}x_{1}-a_{3} & 0\\
-\left(p_{1}x_{2}+b_{1}\right) & -\left(q_{1}x_{1}-a_{1}\right) & 0 & -c_{2} & p_{2}x_{4}+b_{2}\\
-\left(b_{3}+p_{3}x_{2}\right) & -\left(q_{3}x_{1}-a_{3}\right) & c_{2} & 0 & q_{2}x_{3}-a_{2}\\
0 & 0 & -\left(p_{2}x_{4}+b_{2}\right) & -\left(q_{2}x_{3}-a_{2}\right) & 0
\end{array}\right],
\]
and the Jacobi condition simplifies to 
\begin{gather}
\epsilon_{ijk}\mathrm{J}_{im}\frac{\partial\mathrm{J}_{jk}}{\partial x_{m}}=2\left(p_{1}-p_{2}\right)\left(a_{3}-q_{3}x_{1}\right)-2\left(p_{3}+q_{2}\right)\left(a_{1}-q_{1}x_{1}\right)\nonumber \\
+2\left(p_{2}-q_{1}\right)\left(b_{3}+p_{3}x_{2}\right)+2\left(q_{2}+q_{3}\right)\left(b_{1}+p_{1}x_{2}\right)=0.\label{eq:p32}
\end{gather}
For this to hold for all $\mathbf{x}$ the coefficients of $x_{1}$,
$x_{2}$, and the constant terms must each vanish, giving conditions
solved by $p_{1}=p_{2}=q_{1}$ and $p_{3}=q_{3}=-q_{2}$. The Hamiltonian
model is 
\begin{align}
\dot{x}_{1} & =\left(p_{1}x_{2}x_{3}+b_{1}x_{3}-c_{1}x_{2}\right)+\left[-q_{2}x_{2}x_{4}+b_{3}x_{4}-c_{3}x_{2}\right]\nonumber \\
\dot{x}_{2} & =\left(p_{1}x_{3}x_{1}+c_{1}x_{1}-a_{1}x_{3}\right)+\left[-q_{2}x_{4}x_{1}+c_{3}x_{1}-a_{3}x_{4}\right]\nonumber \\
\dot{x}_{3} & =\left(-2p_{1}x_{1}x_{2}+a_{1}x_{2}-b_{1}x_{1}\right)+\left\{ p_{1}x_{4}x_{5}+b_{2}x_{5}-c_{2}x_{4}\right\} \nonumber \\
\dot{x}_{4} & =\left\{ q_{2}x_{5}x_{3}+c_{2}x_{3}-a_{2}x_{5}\right\} +\left[2q_{2}x_{1}x_{2}+a_{3}x_{2}-b_{3}x_{1}\right]\nonumber \\
\dot{x}_{5} & =\left\{ -\left(p_{1}+q_{2}\right)x_{3}x_{4}+a_{2}x_{4}-b_{2}x_{3}\right\} .\label{eq:m3_ham}
\end{align}
Note that omitting the third gyrostat reduces Eq.~(\ref{eq:p32})
to Eq.~(\ref{eq:jac2}), and the new constraint $p_{1}=p_{2}=q_{1}$
imposed here was absent in Model~2. This cross-coupling between gyrostat
parameters, whereby adding a gyrostat retroactively constrains earlier
ones, is a characteristic of non-nested coupling topologies. It contrasts
with the nested hierarchy results of Section~\ref{sec:nested_h},
where recurrent Jacobi conditions avoid such retroactive constraints
on gyrostats that earlier appeared in simpler members of the hierarchy.

\subsection{Hamiltonian conditions over all representations of $\mathrm{J}$}

\label{sec:3K}

The Hamiltonian conditions derived above, Eq.~(\ref{eq:ham1_cond})
for Model~1, Eq.~(\ref{eq:ham2_cond}) for Model~2, and Eq.~(\ref{eq:p32})
for Model~3, were obtained by evaluating the Jacobi condition on
the canonical Poisson matrix $\mathrm{J}=\sum_{k}\mathrm{J}^{(k)}$
constructed in Section~\ref{sec:ham}. In that construction each
gyrostat's contribution $\mathrm{L}^{(k)}$ places the nonlinear entries
in the third column and row. However, for any single gyrostat with
modes $\{m_{1},m_{2},m_{3}\}$ and parameters $\{p_{k},q_{k},r_{k},a_{k},b_{k},c_{k}\}$,
there are in fact three distinct skew-symmetric matrices that each
recover the correct vector field via $\dot{x}_{i}=\mathrm{J}_{ij}x_{j}$:
\begin{align}
\mathrm{L}_{A}^{(k)} & =\left[\begin{array}{ccc}
0 & -c_{k} & p_{k}x_{m_{2}}+b_{k}\\
c_{k} & 0 & q_{k}x_{m_{1}}-a_{k}\\
-p_{k}x_{m_{2}}-b_{k} & -q_{k}x_{m_{1}}+a_{k} & 0
\end{array}\right],\label{eq:LA}\\[6pt]
\mathrm{L}_{B}^{(k)} & =\left[\begin{array}{ccc}
0 & p_{k}x_{m_{3}}-c_{k} & b_{k}\\
-p_{k}x_{m_{3}}+c_{k} & 0 & -r_{k}x_{m_{1}}-a_{k}\\
-b_{k} & r_{k}x_{m_{1}}+a_{k} & 0
\end{array}\right],\label{eq:LB}\\[6pt]
\mathrm{L}_{C}^{(k)} & =\left[\begin{array}{ccc}
0 & -q_{k}x_{m_{3}}-c_{k} & -r_{k}x_{m_{2}}+b_{k}\\
q_{k}x_{m_{3}}+c_{k} & 0 & -a_{k}\\
r_{k}x_{m_{2}}-b_{k} & a_{k} & 0
\end{array}\right],\label{eq:LC}
\end{align}
where the rows and columns correspond to modes $m_{1}^{(k)},m_{2}^{(k)},m_{3}^{(k)}$
respectively. The three choices place the nonlinear entries in the
third column/row ($\mathrm{L}_{A}$, the canonical choice of Section~\ref{sec:ham}),
the second column/row ($\mathrm{L}_{B}$), or the first column/row
($\mathrm{L}_{C}$). One verifies that each satisfies $\dot{x}_{i}=\mathrm{J}_{ij}x_{j}$
using $p_{k}+q_{k}+r_{k}=0$, and that each individual $\mathrm{L}^{(k)}$
satisfies the Jacobi condition, consistent with the single gyrostat
being unconditionally Hamiltonian.

For a GLOM with $K$ gyrostats, one may choose any combination of
$\mathrm{L}_{A}^{(k)}$, $\mathrm{L}_{B}^{(k)}$, $\mathrm{L}_{C}^{(k)}$
independently for each $k$, yielding $3^{K}$ candidate matrices
$\mathrm{J}=\sum_{k=1}^{K}\mathrm{J}^{(k)}$. The conditions for the
Jacobi identity $\epsilon_{ijk}\mathrm{J}_{im}\partial\mathrm{J}_{jk}/\partial x_{m}=0$
vary across these $3^{K}$ candidates, and for the GLOM with free
parameters each condition gives rise to a corresponding Hamiltonian
reduction, for a maximum of up to $3^{K}$ sets of conditions since
these conditions might not all be distinct.

\subsubsection*{Results for Models 1 and 2}

For Model~1 ($K=2$) there are $3^{2}=9$ candidate matrices. We
have evaluated the Jacobi condition for each candidate symbolically
(MATLAB code in SI, corrected as described therein). In all nine cases
the condition $\epsilon_{ijk}\mathrm{J}_{im}\partial\mathrm{J}_{jk}/\partial x_{m}=0$
requires the product of two nonlinear parameters (one from each gyrostat)
to vanish; the secondary linear condition ($p_{2}b_{1}-p_{1}b_{2}-q_{2}c_{1}=0$
in the canonical case) has analogues in each representation, all equivalent
under relabelling. In every case the same physical content holds:
at least one gyrostat must have its leading nonlinear term absent.

For Model~2 ($K=2$, sparse coupling), the exhaustive check reveals
a qualitative difference. Eight of the nine representations require
the vanishing of a nonlinear coefficient analogous to $q_{2}=0$ in
the canonical form. One representation, however, is unconditionally
Hamiltonian:

\begin{proposition}[Unconditional Hamiltonian structure of the $(L_{A},L_{C})$
representation]\label{prop:AC_uncond} For the $K=2$ sparse Model~2
with gyrostat~1 on modes $\{x_{1},x_{2},x_{3}\}$ and gyrostat~2
on modes $\{x_{3},x_{4},x_{5}\}$, the Poisson matrix $\mathrm{J}_{(A,C)}$
comprised of $\mathrm{L}_{A}^{(1)}$ and $\mathrm{L}_{C}^{(2)}$ satisfies
the Jacobi identity for all parameter values. The corresponding Casimir
is 
\begin{equation}
C_{(A,C)}=a_{2}\!\left(a_{1}x_{1}-\tfrac{q_{1}}{2}x_{1}^{2}\right)+a_{2}\!\left(b_{1}x_{2}+\tfrac{p_{1}}{2}x_{2}^{2}\right)+a_{2}c_{1}x_{3}+c_{1}\!\left(b_{2}x_{4}-\tfrac{r_{2}}{2}x_{4}^{2}\right)+c_{1}\!\left(c_{2}x_{5}+\tfrac{q_{2}}{2}x_{5}^{2}\right),\label{eq:cas_AC}
\end{equation}
which is valid for all parameter values. No representation of any
$K\geq3$ sparse GLOM is unconditionally Hamiltonian.\end{proposition}
\begin{proof}
In $\mathrm{L}_{C}^{(2)}$ the state-dependent entries involve only
the private modes $\{x_{4},x_{5}\}$ of gyrostat~2 (the entry $\mathrm{J}_{34}^{(2)}=-c_{2}-q_{2}x_{5}$
depends on $x_{5}$ and $\mathrm{J}_{35}^{(2)}=b_{2}-r_{2}x_{4}$
depends on $x_{4}$). In $\mathrm{L}_{A}^{(1)}$ the state-dependent
entries involve only the private modes $\{x_{1},x_{2}\}$ of gyrostat~1.
The two private-mode sets are disjoint: $\{x_{1},x_{2}\}\cap\{x_{4},x_{5}\}=\emptyset$.
The cross-term $\epsilon_{ijk}\mathrm{J}_{im}^{(1)}\partial\mathrm{J}_{jk}^{(2)}/\partial x_{m}$
requires $m\in\{x_{4},x_{5}\}$ (where $\mathrm{J}^{(2)}$ has state
dependence) but $\mathrm{J}_{im}^{(1)}=0$ for $m\in\{x_{4},x_{5}\}$;
the reverse cross-term requires $m\in\{x_{1},x_{2}\}$ but $\mathrm{J}_{im}^{(2)}=0$
for $m\in\{x_{1},x_{2}\}$. Both cross-terms therefore vanish for
every $\mathbf{x}$ and every parameter choice. That $C_{(A,C)}$
is the Casimir follows by computing $\mathrm{NULL}(\mathrm{J}_{(A,C)})$
and integrating the gradient (verified symbolically). The failure
at $K\geq3$ follows because every $L_{C}$ representation of any
gyrostat $k\geq2$ in the sparse hierarchy places a state-dependent
entry involving the shared mode $x_{2k+1}$, which is also a private
mode of gyrostat $k+1$; the cross-terms between adjacent gyrostats
therefore no longer vanish identically for any combination of representations
(as is verified using exhaustive enumeration for $K=3$ in SI).
\end{proof}
Equation~(\ref{eq:cas_AC}) is the first Casimir for a $K=2$ GLOM
with all nonlinear parameters generically nonzero. Comparing with
the canonical Casimir $C_{a}$ (Eq.~\ref{eq:cas_m2}), the two expressions
agree on the gyrostat-1 block $\{x_{1},x_{2},x_{3}\}$ and differ
only in the gyrostat-2 private modes: $+p_{2}x_{4}^{2}/2$ is replaced
by $-r_{2}x_{4}^{2}/2$, and the linear $c_{2}x_{5}$ is replaced
by $c_{2}x_{5}+q_{2}x_{5}^{2}/2$. When $q_{2}=0$ and $p_{2}=-r_{2}$
(the standard Hamiltonian condition) the two Casimirs coincide up
to a sign in the $x_{4}$ term.

\subsubsection*{A contrasting example: EC-LOM(4)}

The necessity of searching all $3^{K}$ candidates is further illustrated
by the well-known energy-conserving model of \citet{Lorenz1996} involving
slow modes on a circle, whose conservative core can be written as
a coupled gyrostat system \citep{Hu2026} denoted EC-LOM($n,m$).
Here the GLOM parameters are fixed by the physical model and the question
is purely whether Hamiltonian structure exists. Taking EC-LOM($4,0$)
with $K=4$ gyrostats, there are $3^{4}=81$ candidate matrices to
test. The canonical choice of $\mathrm{J}$ fails the Jacobi condition
(the residual is $4x_{1}-2x_{2}+2x_{3}\neq0$). An exhaustive symbolic
search over all 81 candidates finds that none satisfies the Jacobi
condition. The residuals for all 81 are nonzero polynomials in $\mathbf{x}$.
The EC-LOM($4,0$) therefore does not admit non-canonical Hamiltonian
structure within the class of gyrostat-sum Poisson matrices; whether
it is Hamiltonian via a matrix not of this form remains an open question.
This negative result itself demonstrates the value of the $3^{K}$
search: without it, the failure of the canonical representation alone
could not rule out Hamiltonian structure in the broader class.

\subsection{Nested Hamiltonian hierarchies}

\label{sec:nested_h}

\subsubsection*{Sparse hierarchy}

We construct hierarchies extending Model~2 for $K=1,2,3,\ldots$
with $M=2K+1$. For $K=1$ the single gyrostat is always Hamiltonian
(Section~\ref{sec:ham}). The incremental Jacobi condition for each
new $K$, obtained by subtracting the condition for $K-1$ from that
for $K$, is 
\begin{align}
K=2: & \quad q_{2}\!\left(q_{1}x_{1}+p_{1}x_{2}+b_{1}-a_{1}\right)=0,\nonumber \\
K=3: & \quad q_{3}\!\left(q_{2}x_{3}+p_{2}x_{4}+b_{2}-a_{2}\right)=0,\nonumber \\
K=4: & \quad q_{4}\!\left(q_{3}x_{5}+p_{3}x_{6}+b_{3}-a_{3}\right)=0,\label{eq:34}
\end{align}
as derived in SI. The simplest Hamiltonian hierarchy takes $q_{k}=0$
for all $k\geq2$, giving at each level $K$ a model of the form 
\begin{align*}
\dot{x}_{1} & =\left(p_{1}x_{2}x_{3}+b_{1}x_{3}-c_{1}x_{2}\right)\\
\dot{x}_{2} & =\left(q_{1}x_{3}x_{1}+c_{1}x_{1}-a_{1}x_{3}\right)\\
\dot{x}_{3} & =\left(r_{1}x_{1}x_{2}+a_{1}x_{2}-b_{1}x_{1}\right)+\left\{ p_{2}x_{4}x_{5}+b_{2}x_{5}-c_{2}x_{4}\right\} \\
\dot{x}_{4} & =\left\{ c_{2}x_{3}-a_{2}x_{5}\right\} \\
\dot{x}_{5} & =\left\{ -p_{2}x_{3}x_{4}+a_{2}x_{4}-b_{2}x_{3}\right\} ,
\end{align*}
and so on, with each new gyrostat $k\geq2$ having nonzero $p_{k},r_{k}=-p_{k}$
and $q_{k}=0$. Of course, this characterization comes from the default
representation of $\mathrm{J}$, and other Hamiltonian hierarchies
are possible upon searching through the $3^{K}$ possibilities (see
Section 4.2).

Table~\ref{tab:casimir_null} lists $\mathrm{NULL}\left(\mathrm{J}\right)$
for $K=1,\ldots,4$ in this hierarchy: every entry is a gradient vector,
and each level acquires has exactly one Casimir. Furthermore, the
Casimir gradients are consistent under projection: restricting the
$K=3$ gradient to the $M=5$ subspace gives 
\[
a_{3}\left[\begin{array}{c}
a_{2}\left(a_{1}-q_{1}x_{1}\right)\\
a_{2}\left(b_{1}+p_{1}x_{2}\right)\\
a_{2}c_{1}\\
c_{1}\left(b_{2}+p_{2}x_{4}\right)\\
c_{1}c_{2}
\end{array}\right],
\]
which is collinear with $\nabla C_{a}$ for $K=2$.

\subsubsection*{Dense hierarchy}

Extending Model~1 to $K=3$ gives, for example, 
\begin{align}
\dot{x}_{1} & =\left(p_{1}x_{2}x_{3}+b_{1}x_{3}-c_{1}x_{2}\right)\nonumber \\
\dot{x}_{2} & =\left(q_{1}x_{3}x_{1}+c_{1}x_{1}-a_{1}x_{3}\right)+\left\{ p_{2}x_{3}x_{4}+b_{2}x_{4}-c_{2}x_{3}\right\} \nonumber \\
\dot{x}_{3} & =\left(r_{1}x_{1}x_{2}+a_{1}x_{2}-b_{1}x_{1}\right)+\left\{ q_{2}x_{4}x_{2}+c_{2}x_{2}-a_{2}x_{4}\right\} +\left[p_{3}x_{4}x_{5}+b_{3}x_{5}-c_{3}x_{4}\right]\nonumber \\
\dot{x}_{4} & =\left\{ r_{2}x_{2}x_{3}+a_{2}x_{3}-b_{2}x_{2}\right\} +\left[q_{3}x_{5}x_{3}+c_{3}x_{3}-a_{3}x_{5}\right]\nonumber \\
\dot{x}_{5} & =\left[r_{3}x_{3}x_{4}+a_{3}x_{4}-b_{3}x_{3}\right],\label{eq:dense3}
\end{align}
introducing one new mode per added gyrostat ($M=K+2$). The incremental
Jacobi conditions for this dense hierarchy and default choice of Poisson
matrix are 
\begin{align}
K=2: & \quad p_{1}p_{2}\left(x_{2}-x_{3}\right)+b_{1}p_{2}-b_{2}p_{1}+q_{2}\left(-c_{1}\right)=0,\nonumber \\
K=3: & \quad p_{2}p_{3}\left(x_{3}-x_{4}\right)+b_{2}p_{3}-b_{3}p_{2}+q_{3}\left(q_{1}x_{1}+p_{1}x_{2}+b_{1}-a_{1}-c_{2}\right)=0,\nonumber \\
K=4: & \quad p_{3}p_{4}\left(x_{4}-x_{5}\right)+b_{3}p_{4}-b_{4}p_{3}+q_{4}\left(q_{2}x_{2}+p_{2}x_{3}+b_{2}-a_{2}-c_{3}\right)=0.\label{eq:dense_jac}
\end{align}
Setting $q_{k}=0$ for $k\geq2$ simplifies these conditions to $p_{k-1}p_{k}=0$
and $b_{k-1}p_{k}=b_{k}p_{k-1}$ at each level. Three distinct Hamiltonian
parameter patterns then arise, corresponding to the nonlinear coefficient
sequences being: all-zero ($p_{k}=0$ for all $k$); alternating-nonzero
with even gyrostats vanishing ($p_{2k}=0$, $b_{2k}=0$ for all $k$);
and alternating-nonzero with odd gyrostats vanishing ($p_{2k-1}=0$,
$b_{2k-1}=0$ for all $k$). Table~\ref{tab:casimir_null} lists
$\mathrm{NULL}\left(\mathrm{J}\right)$ for the latter two patterns.
Both show an alternating structure in which odd-numbered $K$ have
a single Casimir and even-numbered $K$ have none, reflecting the
alternating presence of nonlinear entries in $\mathrm{J}$. Casimir
gradients for odd $K$ project consistently to those for smaller odd
$K$ (with $c_{2}=0$), consistent with the more general results on
consistency below.

\medskip{}

\noindent\textbf{Table~4.} Gradient of Casimir ($\mathrm{NULL}\left(\mathrm{J}\right)$)
for sparse and dense Hamiltonian hierarchies. Sparse hierarchy: $q_{2},q_{3},\ldots=0$.
Dense hierarchy~1: $q_{2},q_{3},\ldots=0$, $p_{2},p_{4},\ldots=0$,
$b_{2},b_{4},\ldots=0$. Dense hierarchy~2: $q_{2},q_{3},\ldots=0$,
$p_{1},p_{3},\ldots=0$, $b_{1},b_{3},\ldots=0$. A zero vector indicates
no Casimir at that level.\label{tab:casimir_null}

\begin{tabular}{|c|c|c|c|}
\hline 
$K$  & Sparse  & Dense~1  & Dense~2\tabularnewline
\hline 
\hline 
$1$  & $\left[\begin{array}{c}
a_{1}-q_{1}x_{1}\\
b_{1}+p_{1}x_{2}\\
c_{1}
\end{array}\right]$  & $\left[\begin{array}{c}
a_{1}-q_{1}x_{1}\\
b_{1}+p_{1}x_{2}\\
c_{1}
\end{array}\right]$  & $\left[\begin{array}{c}
a_{1}-q_{1}x_{1}\\
0\\
c_{1}
\end{array}\right]$\tabularnewline
\hline 
$2$  & $\left[\begin{array}{c}
a_{2}\!\left(a_{1}-q_{1}x_{1}\right)\\
a_{2}\!\left(b_{1}+p_{1}x_{2}\right)\\
a_{2}c_{1}\\
c_{1}\!\left(b_{2}+p_{2}x_{4}\right)\\
c_{1}c_{2}
\end{array}\right]$  & $\left[\begin{array}{c}
0\\
0\\
0\\
0
\end{array}\right]$  & $\left[\begin{array}{c}
0\\
0\\
0\\
0
\end{array}\right]$\tabularnewline
\hline 
$3$  & $\left[\begin{array}{c}
a_{2}a_{3}\!\left(a_{1}-q_{1}x_{1}\right)\\
a_{2}a_{3}\!\left(b_{1}+p_{1}x_{2}\right)\\
a_{2}a_{3}c_{1}\\
a_{3}c_{1}\!\left(b_{2}+p_{2}x_{4}\right)\\
a_{3}c_{1}c_{2}\\
c_{1}c_{2}\!\left(b_{3}+p_{3}x_{6}\right)\\
c_{1}c_{2}c_{3}
\end{array}\right]$  & $\left[\begin{array}{c}
a_{3}\!\left(a_{1}+c_{2}-q_{1}x_{1}\right)\\
a_{3}\!\left(b_{1}+p_{1}x_{2}\right)\\
a_{3}c_{1}\\
c_{1}\!\left(b_{3}+p_{3}x_{4}\right)\\
c_{1}\!\left(a_{2}+c_{3}\right)
\end{array}\right]$  & $\left[\begin{array}{c}
a_{3}\!\left(a_{1}+c_{2}-q_{1}x_{1}\right)\\
0\\
a_{3}c_{1}\\
0\\
c_{1}\!\left(a_{2}+c_{3}\right)
\end{array}\right]$\tabularnewline
\hline 
$4$  & $\left[\begin{array}{c}
a_{2}a_{3}a_{4}\!\left(a_{1}-q_{1}x_{1}\right)\\
a_{2}a_{3}a_{4}\!\left(b_{1}+p_{1}x_{2}\right)\\
a_{2}a_{3}a_{4}c_{1}\\
a_{3}a_{4}c_{1}\!\left(b_{2}+p_{2}x_{4}\right)\\
a_{3}a_{4}c_{1}c_{2}\\
a_{4}c_{1}c_{2}\!\left(b_{3}+p_{3}x_{6}\right)\\
a_{4}c_{1}c_{2}c_{3}\\
c_{1}c_{2}c_{3}\!\left(b_{4}+p_{4}x_{8}\right)\\
c_{1}c_{2}c_{3}c_{4}
\end{array}\right]$  & $\left[\begin{array}{c}
0\\
0\\
0\\
0\\
0\\
0
\end{array}\right]$  & $\left[\begin{array}{c}
0\\
0\\
0\\
0\\
0\\
0
\end{array}\right]$\tabularnewline
\hline 
\end{tabular}

\medskip{}

The consistent pattern of Jacobi conditions in both hierarchies is
an example of a general result:

\begin{proposition}[Recurrence of Hamiltonian conditions]\label{prop:recurse}
Consider nested Hamiltonian hierarchies in which $K$ is incremented
by one at each step with a consistent mode-coupling pattern. The incremental
condition for the enlarged system to remain non-canonical Hamiltonian
depends only on the parameters of the newly added gyrostat and those
of immediately preceding gyrostats that share modes with it. \end{proposition} 
\begin{proof}
Since $\mathrm{J}=\sum_{k=1}^{K}\mathrm{J}^{(k)}$, the Jacobi condition
$\epsilon_{ijk}\mathrm{J}_{im}\partial\mathrm{J}_{jk}/\partial x_{m}=0$
for the $K$-gyrostat system can be expanded as a sum over pairs of
gyrostats. The incremental condition (difference between the $K$
and $K-1$ conditions) receives contributions only from terms involving
$\mathrm{J}^{(K)}$: specifically $\epsilon_{ijk}\sum_{l=1}^{K-1}\mathrm{J}_{im}^{(l)}\partial\mathrm{J}_{jk}^{(K)}/\partial x_{m}$
and $\epsilon_{ijk}\sum_{l=1}^{K-1}\mathrm{J}_{im}^{(K)}\partial\mathrm{J}_{jk}^{(l)}/\partial x_{m}$.
Since $\partial\mathrm{J}_{jk}^{(K)}/\partial x_{m}\neq0$ only when
$m$ is one of the modes of gyrostat $K$, only those previous gyrostats
$l$ whose mode sets overlap with gyrostat $K$ contribute to the
sum. For nested hierarchies with consistent coupling, this overlap
set is the same at each $K$, giving the recurrence. The full demonstration
for the sparse and dense cases is in Appendix~3. 
\end{proof}

\subsection{Coupled Hamiltonian hierarchies}

\label{sec:coupled_h}

\subsubsection*{Hub-coupled hierarchy (Model~4)}

In 2D Rayleigh--B�nard convection (Model~4), the three gyrostats
each share mode $x_{1}$, forming a hub-spoke coupling topology with
$K=3$, $M=7$ (Eq.~\ref{eq:p38}). This contrasts with the nested
sparse hierarchy, where each new gyrostat shares one mode with its
immediate predecessor and introduces two fresh modes. Here the additional
conditions for Hamiltonian structure, from the incremental Jacobi
condition, are 
\begin{alignat}{1}
K=2: & \quad{-q_{1}p_{2}x_{4}+q_{1}\left(c_{2}-b_{2}\right)-q_{2}p_{1}x_{2}+\left(c_{1}-b_{1}\right)q_{2}=0,}\nonumber \\
K=3: & \quad{-\left(q_{1}+q_{2}\right)p_{3}x_{6}+\left(q_{1}+q_{2}\right)\left(c_{3}-b_{3}\right)-q_{3}\left(p_{1}x_{2}+p_{2}x_{4}\right)+q_{3}\left(c_{1}-b_{1}+c_{2}-b_{2}\right)=0.}\label{eq:p40}
\end{alignat}
In contrast to the sparse hierarchy (Eq.~\ref{eq:34}), these conditions
are not recurrent: the $K=3$ condition involves parameters of all
three gyrostats. Model~4 itself is not Hamiltonian \citep{Gluhovsky2006},
and its invariants must be found by the standard approach. However,
the specialized model 
\begin{align}
\dot{x}_{1} & =\left(p_{1}x_{2}x_{3}+b_{1}x_{3}-c_{1}x_{2}\right)+\left\{ b_{2}x_{5}-b_{2}x_{4}\right\} +\left[b_{3}x_{7}-b_{3}x_{6}\right]\nonumber \\
\dot{x}_{2} & =\left(q_{1}x_{3}x_{1}+c_{1}x_{1}-a_{1}x_{3}\right)\nonumber \\
\dot{x}_{3} & =\left(r_{1}x_{1}x_{2}+a_{1}x_{2}-b_{1}x_{1}\right)\nonumber \\
\dot{x}_{4} & =\left\{ b_{2}x_{1}-a_{2}x_{5}\right\} \nonumber \\
\dot{x}_{5} & =\left\{ a_{2}x_{4}-b_{2}x_{1}\right\} \nonumber \\
\dot{x}_{6} & =\left[b_{3}x_{1}-a_{3}x_{7}\right]\nonumber \\
\dot{x}_{7} & =\left[a_{3}x_{6}-b_{3}x_{1}\right]\label{eq:p41}
\end{align}
(obtained by setting $p_{2}=q_{2}=0$, $c_{2}=b_{2}$ and $p_{3}=q_{3}=0$,
$c_{3}=b_{3}$) is Hamiltonian for $K=1,2,3$. For $K=1$ the nullspace
yields the standard single-gyrostat Casimir. For $K=2$ and $K=3$
the nullspace vectors are not gradient vectors (verified in SI), so
no Casimir exists. Further specialisation of linear parameters could
recover Casimirs for these levels; the present example illustrates
that fully-coupled topologies require more stringent conditions to
maintain Casimir structure across the hierarchy than nested topologies
do.

\subsubsection*{Fully coupled hierarchy (Model~5)}

The 3D Rayleigh--B�nard model (Eq.~\ref{eq:p42}) involves $K=5$
gyrostats coupling modes through both shared mode $x_{1}$ (gyrostats~1
and~2) and through the cross-coupling mode $x_{7}$ (gyrostats~3,
4, and~5). The incremental Jacobi conditions as each gyrostat is
added are (SI): 
\begin{alignat}{1}
K=2: & \quad-q_{1}p_{2}x_{4}+q_{1}\left(c_{2}-b_{2}\right)-q_{2}p_{1}x_{2}+\left(c_{1}-b_{1}\right)q_{2}=0,\nonumber \\
K=3: & \quad0=0\quad\text{(automatically satisfied),}\nonumber \\
K=4: & \quad p_{2}\left(a_{4}+c_{4}\right)-p_{3}\left(a_{4}-b_{4}\right)-p_{4}\left(a_{2}+c_{2}\right)-q_{4}\left(a_{1}-b_{1}\right)\nonumber \\
 & \quad+\left(p_{4}q_{2}+q_{1}q_{4}\right)x_{1}+p_{1}q_{4}x_{2}-\left(p_{2}-p_{3}\right)q_{4}x_{3}+p_{3}p_{4}x_{4}=0,\nonumber \\
K=5: & \quad p_{1}\left(c_{5}-b_{5}\right)+p_{3}\left(b_{5}-a_{5}\right)+p_{5}\left(a_{2}+b_{2}\right)+q_{5}\left(a_{1}-c_{1}\right)\nonumber \\
 & \quad-\left(q_{1}q_{5}+q_{2}p_{5}\right)x_{1}+p_{3}q_{5}x_{2}+p_{2}p_{5}x_{4}+\left(p_{3}-p_{1}\right)p_{5}x_{5}=0.\label{eq:p43}
\end{alignat}
Unlike the sparse hierarchy, adding gyrostat~4 here introduces retroactive
constraints on gyrostat~2 ($a_{2}=-c_{2}$) and gyrostat~1 ($a_{1}=b_{1}$),
and adding gyrostat~5 constrains gyrostat~2 ($b_{2}=-a_{2}$) and
gyrostat~1 ($c_{1}=a_{1}$). This is the phenomenon already noted
for Model~3: when a new gyrostat couples modes belonging to distinct
earlier gyrostats, it imposes constraints on those gyrostats retroactively.

A Hamiltonian instance of the full hierarchy satisfying Eq.~(\ref{eq:p43})
is obtained by setting 
\begin{equation}
q_{1}=q_{2}=p_{3}=q_{4}=p_{5}=0,\quad a_{4}=b_{4}=-c_{4},\quad a_{2}=-c_{2},\quad a_{1}=b_{1},\quad a_{5}=b_{5}=c_{5},\quad b_{2}=-a_{2},\quad c_{1}=a_{1},\label{eq:p44}
\end{equation}
giving the model 
\begin{align}
\dot{x}_{1} & =\left(p_{1}x_{2}x_{3}+a_{1}x_{3}-a_{1}x_{2}\right)+\left\{ p_{2}x_{4}x_{5}-a_{2}x_{5}+a_{2}x_{4}\right\} \nonumber \\
\dot{x}_{2} & =\left(a_{1}x_{1}-a_{1}x_{3}\right)+\left|a_{5}x_{7}-a_{5}x_{5}\right|\nonumber \\
\dot{x}_{3} & =\left(-p_{1}x_{1}x_{2}+a_{1}x_{2}-a_{1}x_{1}\right)+\left\llbracket p_{4}x_{4}x_{7}+a_{4}x_{7}+a_{4}x_{4}\right\rrbracket \nonumber \\
\dot{x}_{4} & =\left\{ -a_{2}x_{1}-a_{2}x_{5}\right\} +\left\llbracket -a_{4}x_{3}-a_{4}x_{7}\right\rrbracket \nonumber \\
\dot{x}_{5} & =\left\{ -p_{2}x_{1}x_{4}+a_{2}x_{4}+a_{2}x_{1}\right\} +\left|q_{5}x_{7}x_{2}+a_{5}x_{2}-a_{5}x_{7}\right|\nonumber \\
\dot{x}_{6} & =\left[b_{3}x_{8}-c_{3}x_{7}\right]\nonumber \\
\dot{x}_{7} & =\left[q_{3}x_{8}x_{6}+c_{3}x_{6}-a_{3}x_{8}\right]+\left\llbracket -p_{4}x_{3}x_{4}+a_{4}x_{4}-a_{4}x_{3}\right\rrbracket +\left|-q_{5}x_{2}x_{5}+a_{5}x_{5}-a_{5}x_{2}\right|\nonumber \\
\dot{x}_{8} & =\left[-q_{3}x_{6}x_{7}+a_{3}x_{7}-b_{3}x_{6}\right],\label{eq:p45}
\end{align}
whose non-canonical Hamiltonian form has been verified in SI. Table~\ref{tab:cas5}
lists $\mathrm{NULL}\left(\mathrm{J}\right)$ for each sub-hierarchy.
Gyrostats~1 and~2 introduce new modes, so their Casimir gradients
are consistent under projection as in the nested case. Gyrostat~3
couples three new modes and introduces a second independent Casimir.
Gyrostats~4 and~5 couple only existing modes; adding them eliminates
all Casimirs, because the cross-coupling between existing modes exhausts
the nullspace of $\mathrm{J}$.

\medskip{}

\noindent\textbf{Table~5.} Gradient of Casimir ($\mathrm{NULL}\left(\mathrm{J}\right)$)
for the fully coupled Hamiltonian hierarchy of Eq.~(\ref{eq:p45}).
A zero vector indicates no Casimir.\label{tab:cas5}

\begin{tabular}{|c|c|}
\hline 
$K$  & $\mathrm{NULL}\left(\mathrm{J}\right)$\tabularnewline
\hline 
\hline 
$1$  & $\left[\begin{array}{c}
a_{1}\\
a_{1}+p_{1}x_{2}\\
a_{1}
\end{array}\right]$\tabularnewline
\hline 
$2$  & $\left[\begin{array}{c}
a_{1}a_{2}\\
a_{2}\left(a_{1}+p_{1}x_{2}\right)\\
a_{1}a_{2}\\
-a_{1}\left(a_{2}-p_{2}x_{4}\right)\\
-a_{1}a_{2}
\end{array}\right]$\tabularnewline
\hline 
$3$  & $\left[\begin{array}{c}
a_{1}a_{2}\\
a_{2}\left(a_{1}+p_{1}x_{2}\right)\\
a_{1}a_{2}\\
-a_{1}\left(a_{2}-p_{2}x_{4}\right)\\
-a_{1}a_{2}\\
0\\
0\\
0
\end{array}\right]$,\quad{}$\left[\begin{array}{c}
0\\
0\\
0\\
0\\
0\\
a_{3}-q_{3}x_{6}\\
b_{3}\\
c_{3}
\end{array}\right]$\tabularnewline
\hline 
$4$  & $\left[\begin{array}{c}
0\\
0\\
0\\
0\\
0\\
0\\
0\\
0
\end{array}\right]$\tabularnewline
\hline 
$5$  & $\left[\begin{array}{c}
0\\
0\\
0\\
0\\
0\\
0\\
0\\
0
\end{array}\right]$\tabularnewline
\hline 
\end{tabular}

\medskip{}

The behaviour shown in Table~\ref{tab:cas5} reflects a general principle
that is formalised in the following result:

\begin{theorem}[Consistency of Casimir gradients]\label{prop:consist}
Consider non-canonical Hamiltonian hierarchies obtained by incrementing
$K$ by one. Where a Casimir is maintained as $K$ increases, its
gradient is consistent under projection to the subspace of the smaller
model. \end{theorem} 
\begin{proof}
Write $\mathrm{J}_{2}=\mathrm{J}_{1}'+\Delta\mathrm{J}$, where $\mathrm{J}_{1}'\in\mathbb{R}^{m_{2}\times m_{2}}$
embeds the $K=k_{1}$ Poisson matrix with zeros in the new mode rows/columns,
and $\Delta\mathrm{J}$ carries the contributions of the additional
gyrostats. From the well-known result on null-spaces, $\mathrm{NULL}\left(\mathrm{J}_{1}'\right)\cap\mathrm{NULL}\left(\Delta\mathrm{J}\right)\subseteq\mathrm{NULL}\left(\mathrm{J}_{2}\right)$.
Restricting to subspaces that are also gradient vectors (denoted $\mathrm{NULL}^{g}$)
and using the fact that a sum of gradient vectors is a gradient, $\mathrm{NULL}^{g}\left(\mathrm{J}_{1}'\right)\cap\mathrm{NULL}^{g}\left(\Delta\mathrm{J}\right)\subseteq\mathrm{NULL}^{g}\left(\mathrm{J}_{2}\right)$.
If $\dim\!\left(\mathrm{NULL}^{g}\left(\mathrm{J}_{1}'\right)\right)=1$
and the intersection is nontrivial, then $\mathrm{NULL}^{g}\left(\mathrm{J}_{1}'\right)\subseteq\mathrm{NULL}^{g}\left(\mathrm{J}_{2}\right)$,
so the Casimir of the smaller model is also a Casimir of the larger
one, and its gradient projects consistently. The extension to multiple
Casimirs is analogous. The illustration for the sparse hierarchy is
in Appendix~4. 
\end{proof}
Theorem~\ref{prop:consist} is a key property that highlights the
legibility of Hamiltonian hierarchies, whose invariants can be described
in terms of null-spaces of explicitly constructed matrix sums. This
property guarantees that invariants are not merely repeated at each
new level but are substantially maintained by the larger model when
they persist, making Hamiltonian hierarchies suitable as consistent
families of reduced-order models.

\section{Conclusions}

\label{sec:conclusions}

This paper has developed a systematic theory of quadratic invariants
in coupled gyrostat low-order models (GLOMs), examining both the combinatorial
intractability of the standard algebraic approach and the question
of when a geometric alternative that arises from non-canonical Hamiltonian
structure is available.

\paragraph{Main results.}

The central findings can be stated as follows. Energy is the only
guaranteed invariant for generic GLOMs; all other invariants require
parameter constraints that depend sensitively on how the gyrostats
are coupled as well as which gyrostat parameters are nonzero, and
not merely on the number of gyrostats $K$ or modes $M$. Proposition~\ref{prop:monotone}
makes one relationship underlying this behaviour precise: making parameters
nonzero cannot increase the number of invariants, so the general case
always yields the minimum number of invariants. At the other extreme,
sparse GLOMs without linear feedback achieve the maximum possible
number of $(M+1)/2$ independent quadratic invariants (Theorem~\ref{prop:sparse}),
which was demonstrated from a rank calculation on $K$ independent
constraint equations in $M=2K+1$ unknowns. Between these extremes,
the count of invariants depends mainly on which linear feedback terms
are present in ways that are highly sensitive to fine structure of
how the gyrostats are coupled and which parameters are nonzero. We
showed using Models~1 and~2 how two subclasses with the same number
of nonzero parameters can differ in their number of invariants depending
solely on which linear parameters are nonzero.

\paragraph{Hamiltonian structure and the $3^{K}$ search.}

Non-canonical Hamiltonian structure provides a framework that sidesteps
the combinatorial explosion of the standard approach. When the Jacobi
condition is satisfied, all invariants beyond energy arise as Casimir
functions of the $M\times M$ Poisson matrix $\mathrm{J}$, and are
recoverable from its nullspace without exhaustive enumeration of all
subclasses. For a given coupling configuration between the gyrostats,
it is this avoidance of exhaustive subclass enumeration (and possibly
special treatment in case of degeneracies) that offers the greatest
advantage of the Hamiltonian conditions. For the models studied, the
Hamiltonian conditions constrain primarily the nonlinear rather than
the linear coefficients, and we have shown (Section~\ref{sec:3K})
that these conditions must be checked across all $3^{K}$ representations
of the Poisson matrix $\mathrm{J}$, not merely the default representation.
The $3^{K}$ search yields two qualitatively different kinds of result:
in eight of the nine representations of Model~2, Hamiltonian structure
requires a parameter constraint analogous to $q_{2}=0$; but the ninth
representation $(L_{A},L_{C})$ is unconditionally Hamiltonian, providing
a Casimir (Proposition~\ref{prop:AC_uncond}, Eq.~\ref{eq:cas_AC})
that is valid without any restriction on the nonlinear parameters.
For the EC-LOM($4,0$) of \citet{Lorenz1996} the search finds the
opposite: none of the 81 candidates satisfies the Jacobi condition,
establishing that this model is not Hamiltonian within the gyrostat-sum
class of representations of the Poisson matrix. Both findings would
be invisible to the canonical representation alone.

\paragraph{Hierarchy properties.}

The $3^{K}$ possibilities also give rise to a large family of Hamiltonian
hierarchies that can be constructed for these GLOMs. For Hamiltonian
hierarchies, Proposition~\ref{prop:recurse} and Theorem~\ref{prop:consist}
establish two structural properties that are not evident in general
non-Hamiltonian truncations. The recurrence of incremental Jacobi
conditions (Proposition~\ref{prop:recurse}) implies that maintaining
Hamiltonian structure as gyrostats are progressively added requires
only local conditions involving the new gyrostat and those sharing
modes with it, rather than global re-verification at each level. The
consistency of Casimir gradients under projection (Theorem~\ref{prop:consist})
implies that invariants are genuinely inherited by larger models in
terms of shared structure and gradients: the Casimir of a $K$-gyrostat
model is the restriction of the Casimir of any $K'>K$ model in the
same hierarchy. This property, which follows directly from the additive
structure $\mathrm{J}=\sum_{k}\mathrm{J}^{(k)}$ and elementary properties
of null-spaces of matrix sums, provides a geometric criterion for
constructing physically consistent model hierarchies: Hamiltonian
hierarchies guarantee that truncating a model to fewer modes yields
a system with compatible invariants, while it is not evident that
non-Hamiltonian truncations must always provide such controls.

\paragraph{Nested vs.\ coupled topologies.}

The coupled hierarchies (Models~4 and 5) illustrate that each of
the two consistency properties manifest differently depending on how
gyrostats are coupled. Nested hierarchies, where each new gyrostat
introduces fresh modes in a consistent pattern, exhibit full recurrence
and consistent Casimir projection across all levels. Coupled hierarchies,
where new gyrostats link existing modes without introducing new ones,
can break the recurrence and can eliminate Casimirs as the system
grows, because cross-coupling between existing modes exhausts the
nullspace of $\mathrm{J}$. The canonical representation of the 2D
Rayleigh--B�nard model (Model~4) is not naturally Hamiltonian \citep{Gluhovsky2006},
while that of the 3D model (Model~5) also requires explicit parameter
constraints of Eq.~(\ref{eq:p44}). In both cases, the Hamiltonian
framework provides a systematic route to identifying which restrictions
on the model support Casimir structure and which do not. In contrast,
the standard algebraic approach cannot consistently perform this role
even at these moderate dimensions. Of course, it remains important
to check for all $3^{K}$ representations of each of these models,
across the hierarchies considered here. 

\paragraph{Open problems.}

There are several open questions raised by these results:

The classification of GLOMs by Hamiltonian type (canonical, non-canonical,
non-Hamiltonian) and its relationship to coupling topology is not
yet understood beyond the cases studied here. Proposition~\ref{prop:AC_uncond}
identifies the first instance of a possibly graph-theoretic mechanism
underlying Hamiltonian structure: unconditional Hamiltonian structure
of the $(L_{A},L_{C})$ representation of Model~2 follows from the
private-mode sets of the two gyrostats in that representation being
disjoint. This points to a general conjecture: a GLOM admits an unconditionally
Hamiltonian representation if there exists a choice of $\mathrm{L}_{A}$,
$\mathrm{L}_{B}$, or $\mathrm{L}_{C}$ for each gyrostat such that
the resulting state-dependent mode sets are pairwise disjoint across
gyrostats. Such a condition, if it exists, could be deduced from the
mode-coupling graph without evaluating the Jacobi identity. Of course,
it doesn't appear that such a condition needs to be necessary, since
cancellations between nonzero terms can also yield Hamiltonian structure.
However the possibility of identifying sufficient conditions such
as the ones above while probing this conjecture, which would characterise
one class of Hamiltonian GLOMs purely in terms of coupling topology,
is a natural direction for subsequent work. For small $K$ the exhaustive
$3^{K}$ search is feasible and provides the empirical basis for testing
such a characterisation.

For the non-Hamiltonian class of GLOMs, additional invariants can
still exist, and the standard approach demonstrates this for specific
subclasses of Models~1-3. Characterising these invariants requires
tools beyond the Casimir framework. The sensitivity of the invariant
count to fine parameter structure suggests that an alternate representation
of the problem such as in terms of a coupling graph might be more
promising.

Another interesting set of problems pertains to the asymptotic dynamics
on invariant manifolds defined by the intersections of Casimirs, as
well as the onset and structure of chaos for GLOMs that possess Casimirs
and therefore fewer effective degrees of freedom on the constant energy
surface. Implications of such restrictions of GLOMs derived from Galerkin
projection of physical models to those having non-canonical Hamiltonian
structure and the resulting properties of Casimir invariants for statistical
equilibria and stability properties as well as selective decay in
the dissipative setting are important research questions.

More widely, extending the framework of consistent Casimirs to include
the forced-dissipative situation is necessary to connect these results
to practical questions in geophysical modelling across low order model
hierarchies. In the dissipative case, Casimirs no longer constrain
attractors, but the geometric structure they encode may still influence
transient behaviour, predictability horizons, and the design of data
assimilation schemes for reduced-order models. It would be interesting
to probe whether and how Hamiltonian hierarchies constrained in this
manner can advance modeling and forecasting applications as compared
to ad-hoc truncations.

\section*{Declarations of interest}

The authors have no competing interests to declare.

\section*{Acknowledgments}

The authors are grateful to Frank Kwasniok and Vishal Vasan for helpful
discussions.

\section*{Author contributions}

Ashwin K. Seshadri: Conceptualization, Formal analysis, Investigation,
Methodology, Writing -- original draft; S. Lakshmivarahan: Conceptualization,
Writing -- review \& editing.

\section*{Code availability}

All MATLAB code supporting the computations in this paper is available
at https://github.com/akseshadri/GLOM-Hamiltonian (DOI: https://doi.org/10.5281/zenodo.20112613). 

\section*{Funding}

None.

\section*{Appendix~1: Proof of Theorem~\ref{prop:sparse}}

%% ------------------------------------------------------------------

We establish that for sparse GLOMs without linear feedback (Section~\ref{sec:sparse}),
any quadratic invariant $C_{K}=\tfrac{1}{2}\sum_{i}d_{i}x_{i}^{2}+\sum_{i<j}e_{ij}x_{i}x_{j}+\sum_{i}f_{i}x_{i}$
must have $f_{i}=0$ for all $i$ and $e_{ij}=0$ for all $i\neq j$.
With these terms absent the invariant condition reduces to the $K$
linear equations of Theorem~\ref{prop:sparse}, whose solution space
has dimension $M-K=(M+1)/2$.

\paragraph*{Step~1: linear coefficients vanish.}

For a sparse model without linear feedback each term $f_{i}\dot{x}_{i}$
in $\dot{C}_{K}=0$ contributes $f_{i}\,p_{k}x_{j}x_{k'}$ (where
$\{i,j,k'\}$ is the mode triple of the gyrostat that couples $x_{i}$).
Since the vector field contains no linear terms, every quadratic monomial
$x_{j}x_{k'}$ in $\dot{C}_{K}$ must come exclusively from such a
contribution. Moreover, because any two distinct gyrostats in a sparse
model share at most one mode, each product $x_{j}x_{k'}$ appears
exactly once in the full expression for $\dot{C}_{K}$. The coefficient
of that unique term is $f_{i}$, which must therefore vanish.

\paragraph*{Step~2: mixed quadratic coefficients vanish.}

Let $e_{ij}$ be the coefficient of $x_{i}x_{j}$ in $C_{K}$. We
distinguish two cases according to whether $i$ and $j$ belong to
the same gyrostat.

\textit{Case~A: modes $i$ and $j$ are coupled by a gyrostat.} There
exists a gyrostat with mode triple $\{i,j,k\}$. The term $e_{ij}x_{i}\dot{x}_{j}$
is proportional to $e_{ij}x_{i}^{2}x_{k}$, and the term $e_{ij}\dot{x}_{i}x_{j}$
is proportional to $e_{ij}x_{j}^{2}x_{k}$. Both cubic monomials $x_{i}^{2}x_{k}$
and $x_{j}^{2}x_{k}$ involve the third mode $x_{k}$ raised to the
first power, and since in a sparse model no other gyrostat couples
$\{i,k\}$ or $\{j,k\}$, each monomial appears exactly once in $\dot{C}_{K}$.
Setting their coefficients to zero forces $e_{ij}=0$. There are $K\binom{3}{2}=3K$
such pairs $(i,j)$.

\textit{Case~B: modes $i$ and $j$ are not coupled by any gyrostat.}
The term $e_{ij}x_{i}\dot{x}_{j}$ contributes a cubic monomial $e_{ij}x_{i}x_{k'}x_{l'}$,
where $\{j,k',l'\}$ is the mode triple of the gyrostat coupling $x_{j}$.
Since $i$ does not belong to that gyrostat triple, the monomial $x_{i}x_{k'}x_{l'}$
(with $i\notin\{k',l'\}$) can arise from at most two sources in $\dot{C}_{K}$:
from $e_{ik'}$ and from $e_{il'}$. Hence each such monomial generates
an independent constraint on the $e_{ij}$ with $i,j$ not coupled
through a gyrostat. The total number of such monomials is at least
\[
\tfrac{1}{2}\bigl[2(K-1)(M-5)+(M-K+1)(M-3)\bigr]=(3K-2)(K-1),
\]
which for $K\geq2$ exceeds the number of unknown coefficients $\binom{M}{2}-3K=2K(K-1)$.
Since each coefficient $e_{ij}$ induces linearly independent terms
in $\dot{C}_{K}$, the only solution is $e_{ij}=0$.

With $f_{i}=e_{ij}=0$, the condition $\dot{C}_{K}=0$ reduces to
the $K$ equations $p_{k}d_{m_{1}^{(k)}}+q_{k}d_{m_{2}^{(k)}}+r_{k}d_{m_{3}^{(k)}}=0$,
$k=1,\ldots,K$, in $M=2K+1$ unknowns. The coefficient matrix has
full rank for all $K\geq1$ by the diagonal-submatrix argument in
the proof of Theorem~\ref{prop:sparse}: the $K\times K$ submatrix
with columns $\{d_{2},d_{4},\ldots,d_{2K}\}$ is diagonal with nonzero
entries $q_{k}$, so the full constraint matrix has rank $K$, giving
$\dim\!\left(\mathrm{NULL}\right)=M-K=K+1=(M+1)/2$.

\begin{remark} An alternative route to Step~2 uses the sign-flip
symmetries of sparse models without linear feedback. For each gyrostat
with mode triple $\{i,j,k\}$, the transformation $x_{i}\to-x_{i}$,
$x_{j}\to-x_{j}$, $x_{k}\to x_{k}$ is a symmetry of the vector field.
By Appendix~2, any invariant must share this symmetry. A mixed term
$e_{mn}x_{m}x_{n}$ changes sign under the transformation if exactly
one of $m,n$ equals $i$ or $j$, forcing $e_{mn}=0$. The full set
of $2^{K}$ sign-flip symmetries (one sign combination per gyrostat)
forces all $e_{ij}=0$ and all $f_{i}=0$. \end{remark}

%% ------------------------------------------------------------------

\section*{Appendix~2: Symmetries and invariants in GLOMs}

%% ------------------------------------------------------------------

Let $\dot{\mathbf{x}}=\mathbf{f}(\mathbf{x})$ with $\mathbf{x}\in\mathbb{R}^{M}$
and $\mathbf{f}:\mathbb{R}^{M}\to\mathbb{R}^{M}$. A map $\mathcal{S}:\mathbb{R}^{M}\to\mathbb{R}^{M}$
is a \emph{symmetry} of the system if $\mathbf{y}=\mathcal{S}\mathbf{x}$
satisfies the same equations $\dot{\mathbf{y}}=\mathbf{f}(\mathbf{y})$
whenever $\dot{\mathbf{x}}=\mathbf{f}(\mathbf{x})$. For an invertible
$\mathcal{S}$ this is equivalent to 
\begin{equation}
\mathbf{f}(\mathbf{x})=\mathcal{S}^{-1}\mathbf{f}(\mathcal{S}\mathbf{x})\quad\text{for all }\mathbf{x}.\label{eq:sym_cond}
\end{equation}

We restrict to \emph{sign-flip symmetries}: diagonal matrices $\mathcal{S}=\mathrm{diag}(s_{1},\ldots,s_{M})$
with $s_{i}=\pm1$. Since $s_{i}^{2}=1$, we have $\mathcal{S}^{-1}=\mathcal{S}$,
and Eq.~(\ref{eq:sym_cond}) becomes 
\begin{equation}
\dot{x}_{i}=f_{i}(\mathbf{x})=s_{i}f_{i}(\mathcal{S}\mathbf{x}),\quad i=1,\ldots,M.\label{eq:sym_comp}
\end{equation}

\begin{proposition}[Symmetry inheritance]\label{prop:sym_inherit}
If $\mathcal{S}$ is a sign-flip symmetry of $\dot{\mathbf{x}}=\mathbf{f}(\mathbf{x})$
and $C(\mathbf{x})$ is an invariant, then $C(\mathcal{S}\mathbf{x})$
is also an invariant. \end{proposition} 
\begin{proof}
Since $\dot{C}(\mathbf{x})=\sum_{i}\frac{\partial C}{\partial x_{i}}f_{i}(\mathbf{x})=0$
and $s_{i}^{2}=1$, we compute 
\begin{align}
\dot{C}(\mathbf{x}) & =\sum_{i=1}^{M}\frac{\partial C}{\partial x_{i}}f_{i}(\mathbf{x})=\sum_{i=1}^{M}\frac{\partial C}{\partial x_{i}}\,s_{i}f_{i}(\mathcal{S}\mathbf{x})\nonumber \\
 & =\sum_{i=1}^{M}\frac{\partial C}{\partial x_{i}}\,\frac{1}{s_{i}}f_{i}(\mathcal{S}\mathbf{x})=\sum_{i=1}^{M}\frac{\partial C}{\partial(s_{i}x_{i})}f_{i}(\mathcal{S}\mathbf{x})=\frac{d}{dt}C(\mathcal{S}\mathbf{x})\Big|_{\mathbf{x}}.
\end{align}
Hence $\dot{C}(\mathbf{x})=0$ implies $\frac{d}{dt}C(\mathcal{S}\mathbf{x})=0$,
so $C(\mathcal{S}\mathbf{x})$ is conserved. 
\end{proof}
\noindent\textit{Consequence for invariant structure.} If $C(\mathcal{S}\mathbf{x})=C(\mathbf{x})$
is required for every symmetry $\mathcal{S}$ of the system, then
any monomial $x_{i_{1}}\cdots x_{i_{r}}$ that changes sign under
some $\mathcal{S}$ cannot appear in $C$. In particular, a linear
term $f_{i}x_{i}$ changes sign under the transformation $x_{i}\to-x_{i}$
whenever that is a symmetry, forcing $f_{i}=0$. A mixed quadratic
term $e_{ij}x_{i}x_{j}$ changes sign under $x_{i}\to-x_{i}$, $x_{j}\to x_{j}$
(or the reverse) whenever that is a symmetry, forcing $e_{ij}=0$.
The Euler gyrostat and the sparse no-feedback models are invariant
under all sign combinations that flip an even number of modes in each
gyrostat triple, providing enough symmetries to eliminate all linear
and mixed quadratic terms from any invariant. This is consistent with
the results of Section~\ref{sec:standard} and Appendix~1, which
establish the same conclusion by the constraint-equation route.

%% ------------------------------------------------------------------

\section*{Appendix~3: Proof of Proposition~\ref{prop:recurse}}

%% ------------------------------------------------------------------

We derive the incremental Jacobi condition and show that it depends
only on the parameters of gyrostat $K+1$ and those previous gyrostats
that share modes with it.

\paragraph*{Incremental condition.}

Since $\mathrm{J}=\sum_{l=1}^{K}\mathrm{J}^{(l)}$, the Jacobi condition
for the $K$-gyrostat system is 
\begin{equation}
\mathcal{J}_{K}\equiv\epsilon_{ijk}\sum_{l=1}^{K}\mathrm{J}_{im}^{(l)}\sum_{l'=1}^{K}\frac{\partial\mathrm{J}_{jk}^{(l')}}{\partial x_{m}}=0,\label{eq:jacK}
\end{equation}
where the inner sums use distinct dummy indices $l$ and $l'$ to
avoid ambiguity. The incremental condition upon adding gyrostat $K+1$
is $\mathcal{J}_{K+1}-\mathcal{J}_{K}=0$. Expanding, 
\begin{align}
\mathcal{J}_{K+1}-\mathcal{J}_{K} & =\epsilon_{ijk}\Biggl(\sum_{l=1}^{K+1}\mathrm{J}_{im}^{(l)}\sum_{l'=1}^{K+1}\frac{\partial\mathrm{J}_{jk}^{(l')}}{\partial x_{m}}-\sum_{l=1}^{K}\mathrm{J}_{im}^{(l)}\sum_{l'=1}^{K}\frac{\partial\mathrm{J}_{jk}^{(l')}}{\partial x_{m}}\Biggr)\nonumber \\
 & =\epsilon_{ijk}\Biggl(\mathrm{J}_{im}^{(K+1)}\frac{\partial\mathrm{J}_{jk}^{(K+1)}}{\partial x_{m}}+\sum_{l=1}^{K}\mathrm{J}_{im}^{(l)}\frac{\partial\mathrm{J}_{jk}^{(K+1)}}{\partial x_{m}}+\sum_{l'=1}^{K}\mathrm{J}_{im}^{(K+1)}\frac{\partial\mathrm{J}_{jk}^{(l')}}{\partial x_{m}}\Biggr)=0.\label{eq:incr_jac}
\end{align}
The first term in Eq.~(\ref{eq:incr_jac}) is the self-Jacobi condition
of gyrostat $K+1$ alone, which equals zero because every single gyrostat
is non-canonical Hamiltonian (Section~\ref{sec:ham}). The incremental
condition therefore reduces to 
\begin{equation}
\epsilon_{ijk}\sum_{l=1}^{K}\mathrm{J}_{im}^{(l)}\frac{\partial\mathrm{J}_{jk}^{(K+1)}}{\partial x_{m}}+\epsilon_{ijk}\sum_{l'=1}^{K}\mathrm{J}_{im}^{(K+1)}\frac{\partial\mathrm{J}_{jk}^{(l')}}{\partial x_{m}}=0.\label{eq:incr_simp}
\end{equation}

\paragraph*{Locality of the condition.}

Consider the first sum in Eq.~(\ref{eq:incr_simp}). The derivative
$\partial\mathrm{J}_{jk}^{(K+1)}/\partial x_{m}$ is nonzero only
when $m$ equals one of the first two mode indices of gyrostat $K+1$,
namely $m_{1}^{(K+1)}$ or $m_{2}^{(K+1)}$: specifically, 
\begin{alignat}{2}
m=m_{1}^{(K+1)}: & \quad\frac{\partial\mathrm{J}_{jk}^{(K+1)}}{\partial x_{m}} & =\pm q_{K+1} & \quad\text{(from the }q_{K+1}x_{m_{1}}\text{ entries),}\\
m=m_{2}^{(K+1)}: & \quad\frac{\partial\mathrm{J}_{jk}^{(K+1)}}{\partial x_{m}} & =\pm p_{K+1} & \quad\text{(from the }p_{K+1}x_{m_{2}}\text{ entries).}
\end{alignat}
The factor $\mathrm{J}_{im}^{(l)}$ with $m\in\{m_{1}^{(K+1)},m_{2}^{(K+1)}\}$
is nonzero only if gyrostat $l$ involves mode $m_{1}^{(K+1)}$ or
$m_{2}^{(K+1)}$, i.e.\ only if gyrostat $l$ shares a mode with
gyrostat $K+1$. An identical argument applies to the second sum (with
the roles of $l$ and $K+1$ interchanged). Therefore both sums in
Eq.~(\ref{eq:incr_simp}) receive contributions only from those previous
gyrostats $l\in\{1,\ldots,K\}$ whose mode sets overlap with that
of gyrostat $K+1$.

\paragraph*{Recurrence.}

For nested hierarchies with a consistent coupling pattern, the set
of previous gyrostats sharing modes with gyrostat $K+1$ is the same
function of $K$ at every level. In the sparse hierarchy, gyrostat
$K+1$ shares only mode $x_{2K+1}$ with gyrostat $K$ (and no modes
with gyrostats $1,\ldots,K-1$), so only the parameters of gyrostat
$K$ appear in the incremental condition --- giving the single-step
recurrence of Eq.~(\ref{eq:34}). In the dense hierarchy, gyrostat
$K+1$ shares modes $x_{K+1}$ and $x_{K+2}$ with gyrostats $K$
and $K-1$, so parameters of both enter the condition --- giving
the two-step recurrence of Eq.~(\ref{eq:dense_jac}). In both cases
the structure of the incremental condition is the same for every $K$,
establishing Proposition~\ref{prop:recurse}.

%% ------------------------------------------------------------------

\section*{Appendix~4: Illustration of Theorem~\ref{prop:consist}}

%% ------------------------------------------------------------------

We verify Theorem~\ref{prop:consist} concretely for the sparse hierarchy
at the step $k_{1}=1\to k_{2}=2$ (corresponding to $m_{1}=3$ and
$m_{2}=5$ modes). The Poisson matrix of the $K=1$ member, embedded
in $\mathbb{R}^{5\times5}$ with zeros in the new mode rows and columns,
is 
\[
\mathrm{J}_{1}'=\left[\begin{array}{ccccc}
0 & -c_{1} & b_{1}+p_{1}x_{2} & 0 & 0\\
c_{1} & 0 & q_{1}x_{1}-a_{1} & 0 & 0\\
-b_{1}-p_{1}x_{2} & a_{1}-q_{1}x_{1} & 0 & 0 & 0\\
0 & 0 & 0 & 0 & 0\\
0 & 0 & 0 & 0 & 0
\end{array}\right],
\]
and the increment contributed by gyrostat~2 (with $q_{2}=0$) is
\[
\Delta\mathrm{J}=\left[\begin{array}{ccccc}
0 & 0 & 0 & 0 & 0\\
0 & 0 & 0 & 0 & 0\\
0 & 0 & 0 & -c_{2} & b_{2}+p_{2}x_{4}\\
0 & 0 & c_{2} & 0 & -a_{2}\\
0 & 0 & -b_{2}-p_{2}x_{4} & a_{2} & 0
\end{array}\right].
\]
Their nullspaces are spanned by the columns of 
\[
\mathrm{V}_{1}=\left[\begin{array}{ccc}
a_{1}-q_{1}x_{1} & 0 & 0\\
b_{1}+p_{1}x_{2} & 0 & 0\\
c_{1} & 0 & 0\\
0 & 1 & 0\\
0 & 0 & 1
\end{array}\right],\qquad\mathrm{V}_{d}=\left[\begin{array}{ccc}
1 & 0 & 0\\
0 & 1 & 0\\
0 & 0 & a_{2}\\
0 & 0 & b_{2}+p_{2}x_{4}\\
0 & 0 & c_{2}
\end{array}\right],
\]
respectively. Their intersection $\mathrm{NULL}\!\left(\mathrm{J}_{1}'\right)\cap\mathrm{NULL}\!\left(\Delta\mathrm{J}\right)$
consists of all $\mathbf{v}$ expressible as both $\mathbf{v}=\mathrm{V}_{1}\mathbf{x}$
and $\mathbf{v}=\mathrm{V}_{d}\mathbf{y}$, or equivalently the null-space
of $\mathrm{V}=\bigl[\mathrm{V}_{1}\;\;{-\mathrm{V}_{d}}\bigr]$.
Direct computation gives the unique (up to scalar) null vector 
\[
\begin{bmatrix}\mathbf{x}\\
\mathbf{y}
\end{bmatrix}=\begin{bmatrix}a_{2}\\
c_{1}(b_{2}+p_{2}x_{4})\\
c_{1}c_{2}\\
a_{2}(a_{1}-q_{1}x_{1})\\
a_{2}(b_{1}+p_{1}x_{2})\\
c_{1}
\end{bmatrix},
\]
so the common intersection vector is 
\begin{equation}
\mathrm{V}_{1}\mathbf{x}=\mathrm{V}_{d}\mathbf{y}=\left[\begin{array}{c}
a_{2}(a_{1}-q_{1}x_{1})\\
a_{2}(b_{1}+p_{1}x_{2})\\
a_{2}c_{1}\\
c_{1}(b_{2}+p_{2}x_{4})\\
c_{1}c_{2}
\end{array}\right].\label{eq:A4_result}
\end{equation}
This is precisely $\mathrm{NULL}\!\left(\mathrm{J}_{2}\right)$ as
listed in Table~4, confirming that the $K=2$ Casimir gradient lies
in the intersection of the two embedded nullspaces. Since the intersection
vector is a gradient (its components are partial derivatives of the
function $-\tfrac{1}{2}a_{2}q_{1}x_{1}^{2}+\tfrac{1}{2}a_{2}p_{1}x_{2}^{2}+a_{2}c_{1}x_{3}+\ldots$),
Theorem~\ref{prop:consist} applies and the $K=1$ Casimir gradient
projects consistently onto the $K=2$ Casimir gradient: restricting
Eq.~(\ref{eq:A4_result}) to its first three components and factoring
out $a_{2}$ recovers $\nabla C_{a}$ for $K=1$ exactly.

\pagebreak{}

\bibliographystyle{agufull08}
\bibliography{VG_constants}

\end{document}